\documentclass[a4paper]{article}

\usepackage{amsfonts, amssymb, amsmath, theorem}
\usepackage{fullpage}
\usepackage{mathrsfs, dsfont}

\theoremheaderfont\scshape
\newtheorem{theorem}{Theorem}[section]
\newtheorem{proposition}[theorem]{Proposition}
\newtheorem{assumption}[theorem]{Assumption}
\newtheorem{corollary}[theorem]{Corollary}
\newtheorem{lemma}[theorem]{Lemma}
\theorembodyfont{\rmfamily}
\newtheorem{definition}[theorem]{Definition}
\newtheorem{remark}[theorem]{Remark}

\newtheorem{example}[theorem]{Example}
\newenvironment{proof}[1][]{\noindent\textit{Proof#1.} }{\vskip\baselineskip}

\renewcommand\hat[1]{\widehat{#1}}
\newcommand\A[1]{\mathcal A_T^{#1}}
\newcommand\C{\mathcal C}
\newcommand\D{\mathcal D}
\newcommand\F{\mathscr F}

\newcommand\NN{\mathbb N}
\renewcommand\P{\mathcal P}
\newcommand\PP{\mathbb P}
\newcommand\RR{\mathbb R}
\newcommand\U{\mathfrak U}
\newcommand\UU{\mathbb U}
\newcommand\V{\mathfrak V}
\newcommand\X{\mathscr X}

\newcommand\Z{\mathcal Z}

\newcommand\presucc{\succeq}
\newcommand\prepreq{\preceq}
\newcommand\one{\mathds1}
\newcommand\interior{\operatorname{int}}
\newcommand\ri{\operatorname{ri}}
\newcommand\dom{\operatorname{dom}}
\newcommand\aff{\operatorname{aff}}
\newcommand\set[1]{\left\{#1\right\}} 
\newcommand\sets[2]{\set{#1\,:\,#2}} 
\renewcommand\d{\mathrm{d}}

\newcommand\E[2][]{\mathrm E_{#1}\left[ #2\right]}
\newcommand\bigE[2][]{\operatorname E_{#1}\big[ #2\big]}

\newcommand\ip[2]{\left\langle #1,#2\right\rangle}
\newcommand\bigip[2]{\big\langle #1,#2\big\rangle}
\newcommand\radnik[1][\mathbb Q]{\frac{\d #1}{\d\mathbb P}}
\newcommand\cl{\operatorname{cl}}

\newcommand\hypo{\operatorname{hypo}}
\newcommand\ba{\operatorname{ba}}
\newcommand\ca{\operatorname{ca}}

\newcommand\esssup{\operatorname{ess\;sup}}
\newcommand\proj{P}
\newcommand\qed{\hfill$\Box$}

\numberwithin{equation}{section}

\begin{document}

\title{Multivariate Utility Maximization with\\ Proportional Transaction Costs}
\author{Luciano Campi\thanks{CEREMADE, Universit\'e Paris Dauphine.}\and Mark P. Owen\thanks{Maxwell Institute for Mathematical Sciences, and Department of Actuarial Mathematics and Statistics, Heriot-Watt University.}}
\date{}
\maketitle

\begin{abstract}
We present an optimal investment theorem for a currency exchange model with random and possibly discontinuous proportional transaction costs. The investor's preferences are represented by a multivariate utility function, allowing for simultaneous consumption of any prescribed selection of the currencies at a given terminal date. We prove the existence of an optimal portfolio process under the assumption of asymptotic satiability of the value function. Sufficient conditions for asymptotic satiability of the value function include reasonable asymptotic elasticity of the utility function, or a growth condition on its dual function. We show that the portfolio optimization problem can be reformulated in terms of maximization of a terminal liquidation utility function, and that both problems have a common optimizer.



\textbf{AMS Classification (2000):} Primary -- 91B28, Secondary -- 49N15, 49J40, 49J55.
\end{abstract}

\section{Introduction}

In this paper we consider a portfolio optimization problem over a finite time horizon $[0,T]$ in a continuous-time financial market, where an agent can trade between finitely many risky assets with proportional transaction costs. The underlying financial market model is very general; the terms of each trade are described by a bid-ask process $(\Pi_t)_{t\in[0,T]}$ as in \cite{CampScha}, so that transaction costs can be time-dependent, random and have jumps. In this setting, the portfolio process $(V_t)_{t\in[0,T]}$ is a vector-valued process describing at every instant the number of physical units of each asset held by the agent. The example that the reader should always have in mind is an exchange market with $D$ currencies, in which $V_t=(V_t^1,\dots,V_t^D)$ represents how many dollars, euros, pounds and so on, the agent holds at time $t$. The agent is permitted to dynamically rebalance their portfolio within the set of all admissible self-financing portfolio processes as in \cite{CampScha}. To avoid arbitrage, we assume the existence of a strictly consistent pricing system (SCPS) throughout the paper. Precise details and further assumptions about the modelling of the economy are given in Section \ref{sec:prelim}.

We consider an agent who may consume a prescribed selection of the $D$ assets at time $T$. Without loss of generality, we assume that the agent wishes to consume the first $d$ assets, where $1\le d\le D$. We have two main cases in mind namely $d=D$, whereby the investor can consume all assets, and $d=1$, whereby the investor must liquidate to a reference asset immediately prior to consumption. In the latter case, those assets which are not consumed play the role of pure investment assets. We model the agent's preferences towards terminal consumption by means of a multivariate utility function, $U:\RR^d\rightarrow[-\infty,\infty)$, supported on the non-negative orthant $\RR_+^d$ (see Definition \ref{def:utility}). The utility function is assumed to satisfy the following conditions.
\begin{assumption}\label{ass:basic}
\begin{enumerate}
 \item $U$ is upper semi-continuous;
 \item $U$ is strictly concave on the interior of $\RR_+^d$;
 \item $U$ is essentially smooth, i.e. differentiable in the interior of $\RR_+^d$, and its gradient diverges at the boundary of $\RR_+^d$ (see Definition \ref{def:ess_smooth});
 \item $U$ is asymptotically satiable, i.e. there exist positions in the traded assets for which the marginal utility of $U$ can be made arbitrarily small (see Definition \ref{def:ass_satiable}).
\end{enumerate}
\end{assumption}
In the univariate case ($d=1$) the assumption of both essential smoothness and asymptotic satiability is equivalent to the familiar assumption of continuous differentiability together with the Inada conditions $U'(0)=\infty$ and $U'(\infty)=0$. Precise details about the above conditions can be found within Section \ref{sec:prelim}.

In order to express the investor's preferences towards consumption of the first $d$ assets within the setting of the larger economy we adopt the approach of \cite{Kami01}, extending the utility function $U$ to all $D$ assets. We define $\tilde U:\RR^D\rightarrow[-\infty,\infty)$ by
\begin{equation}\label{eqn:extension}
\tilde U(x):=\begin{cases}U(x_1,\dots,x_d),&x\in\RR_+^D\\-\infty,&\text{otherwise.}\end{cases}
\end{equation}
Although the extended utility function $\tilde U$ theoretically models the possibility of consumption of all $D$ assets, the investor has no incentive to consume anything other than the first $d$ assets because the utility is invariant with respect to increased consumption of the remaining $D-d$ assets.

The investor's primal optimization problem is formulated in terms of the \emph{value function} $u:\RR^D\rightarrow[-\infty,\infty]$ defined\footnote{Since $U$ is assumed to be upper semi-continuous, it is Borel measurable. In fact, the assumption that $U$ is upper semi-continuous can be relaxed to Borel measurability throughout the paper, with the exception of Section \ref{sec:liquidation}. We use the standard convention that $\E{\smash{\tilde U(X)}}=-\infty$ whenever $\E{\smash{\tilde U(X)^-}}= \infty$.}
by
\begin{equation}\label{eqn:primal}
u(x):=\sup\sets{\E{\smash{\tilde U(X)}}}{X\in\A{x}},
\end{equation}
where $x$ represents an initial portfolio, and $\A x$ denotes the set of all terminal values of admissible portfolio processes with initial portfolio $x$. Let $\dom(u):=\sets{x\in\RR^D}{u(x)>-\infty}$ denote the effective domain of $u$, and let $\cl(\dom(u))$ and $\interior(\dom(u))$ denote respectively the closure and interior of the effective domain of $u$. The following assumption holds throughout the paper.

\begin{assumption}\label{ass:proper_concave}
$u(x)<\infty$ for some $x\in\interior(\dom(u))$.
\end{assumption}

Our main results are as follows. In Proposition \ref{thm:value_function} we show that (under Assumption \ref{ass:proper_concave}) the value function is a also utility function. We give an explicit characterisation of $\cl(\dom(u))$ in terms of the cone of deterministic terminal portfolios attainable at zero cost. The set $\cl(\dom(u))$ is itself a closed convex cone which strictly contains $\RR_+^D$, reflecting the rather obvious fact that even with an initial short position in some of the assets, the investor may use other positive initial holdings to trade to a terminal position in which they hold non-negative amounts of each asset. In Proposition \ref{thm:duality} we establish a relationship between the primal problem of utility maximization and an appropriate dual minimization problem \eqref{eqn:dual}. The domain of the dual problem is contained in a space of Euclidean vector measures, in contrast to the frictionless case where real-valued measures suffice. We show that the dual problem has a solution whenever $x\in\interior(\dom(u))$. Finally, in Theorem \ref{thm:optimal}, we prove that the utility maximization problem \eqref{eqn:primal} admits a unique solution for all $x\in\interior(\dom(u))$, under the following assumption.
\begin{assumption}\label{ass:additional}
  $u$ is asymptotically satiable (see Definition \ref{def:ass_satiable}).
\end{assumption}

In Corollary \ref{thm:RAEimplAS} we provide sufficient conditions on the utility function $U$ for Assumption \ref{ass:additional} to hold. Also, to place our optimization problem into the context of other papers which require liquidation of terminal portfolios into a reference asset, we show in in Proposition \ref{thm:liquid-equiv} that the utility maximization problem \eqref{eqn:primal} can be reformulated in terms of maximization of a liquidation utility functional. In Proposition \ref{thm:liquid-attained} we show that both formulations of the optimization problem essentially share a common optimizer.

\vspace\baselineskip

Utility maximization problems in markets with transaction costs have been investigated by many authors, typically using either the dynamic programming approach or the martingale duality approach. While the dynamic programming approach is particularly well suited to treating optimization problems with a Markovian state process (see e.g. \cite{DaviNorm,ShreSone}), the duality approach has the advantage that it is applicable to very general models. The first paper to use the duality approach in the setting of proportional transaction costs was \cite{CvitKara}. Cvitani\'c and Karatzas model two assets (a bond and a stock) as It\^o processes, and assume constant proportional transaction costs. At the close of trading they assume that the investor liquidates their portfolio to the bond in order to consume their wealth. In this setting they prove the existence of a solution to the problem of utility maximization, under the assumption that a dual minimization problem admits a solution. The existence of a solution to the dual problem was subsequently proved in \cite{CvitWang}.

In \cite{Kaba}, a much more general formulation of a transaction costs model for a currency market was introduced, based on the key concept of solvency cone. In the same paper, Kabanov also considers the problem of expected utility maximization, with liquidation of the terminal portfolio to a chosen reference currency, which is used throughout as the num\'eraire. Similarly to \cite{CvitKara}, Kabanov proves the existence of an optimal strategy under the assumption that a dual minimization problem admits a solution.

Developments in the generality of Kabanov's transaction costs model in continuous time have since been given in \cite{KabaLast}, where a square-integrability condition was replaced by an admissibility condition, followed by \cite{KabaStri} which treated the case of time-dependent, random transaction costs, provided the solvency cones can be generated by a countable family of continuous processes. More recently, in \cite{CampScha}, Kabanov's model of currency exchange was further developed to allow discontinuous bid-ask processes, and our optimization problem is set within this very general framework.

A important issue for utility maximization under transaction costs is the consideration of how an investor measures their wealth, and thus their utility. In the frictionless case it is normally assumed that there is a single consumption asset, which is used as a num\'eraire (there are exceptions, e.g. \cite{Lakn}). However, in the transaction cost setting it is quite natural to assume that the investor has access to several non-substitutable consumption assets. This is particularly relevant when one considers a model of currency exchange, where there may be, for example, one consumption asset denominated in each currency. Modeling preferences with respect to several consumption assets clearly requires the use of a multivariate utility function.

In \cite{DeelPhamTouz}, Deelstra et al. investigate a utility maximization problem within the transaction costs framework of \cite{KabaLast}. The agent's preferences are described by a multivariate utility function $U$ which is supported on a constant solvency cone. The utility function is not assumed to be smooth so that liquidation can be included as a particular case. In fact, by assuming that the utility function is supported on the solvency cone, \cite{DeelPhamTouz} are implicitly modeling the occurrence of at least one more trade (e.g. liquidation, or an extended trading period) which takes place either on or after the terminal date, but prior to consumption of wealth.

In \cite{Kami01,Kami04}, Kamizono investigates a utility maximization problem which is also set within the transaction costs framework of \cite{KabaLast}. Kamizono argues convincingly that a distinction should be drawn between \emph{direct} utility (i.e. utility derived explicitly from consumption) and \emph{indirect} utility, which depends on further trading, e.g. liquidation. He argues that \cite{DeelPhamTouz} are using a kind of indirect utility function, which is why they need to consider the case of a non-smooth utility function. We choose to adopt the approach of Kamizono in the current paper by using a direct utility function $U$, which is supported on $\RR_+^d$, in the formulation of the primal problem. The value function $u$, defined in \eqref{eqn:primal}, is then a type of indirect utility, whose support (the closure of its effective domain) is intimately connected to the transaction costs structure, as we shall see in Proposition \ref{thm:value_function}. In Example \ref{exa:nonstrictC1} we demonstrate that the value function $u$ may fail to be either strictly concave or differentiable on $\interior(\dom(u))$.

In order to prove the existence of an optimizer in the multivariate setting, most existing papers make fairly strong technical assumptions on the utility function, which do not admit easy economical interpretations. For example, in \cite{DeelPhamTouz,Kami01,Kami04} the utility function is assumed to be bounded below, and unbounded above. In addition, in \cite{DeelPhamTouz} the dual of the utility function is assumed to explode on the boundary of its effective domain, or to be extendable to a neighbourhood of its original domain. In the current paper, Assumption \ref{ass:basic} is the only assumption we shall make directly on the utility function $U$. It is worth noting that, with the exception of Section \ref{sec:liquidation}, the assumption of upper semi-continuity is only used to ensure that $U$ is Borel measurable, and hence that the primal problem \eqref{eqn:primal} is well defined.

A relatively recent development in the theory of utility maximization is the replacement of the assumption of reasonable asymptotic elasticity on the utility function by a weaker condition. In the frictionless setting, \cite{KramScha03} showed that finiteness of the dual of the value function is sufficient for the existence of an optimal portfolio. Since then \cite{BoucPham} have investigated this further under the discrete time model of transaction costs given in \cite{Scha} and \cite{KabaStriRaso}. They prove the existence of an optimal consumption investment strategy under the assumption of finiteness of the convex dual of the value function corresponding to an auxiliary \emph{univariate} primal problem. The reason why \cite{BoucPham} have to employ an auxiliary, univariate primal problem is that the generalization of the methods of Kramkov and Schachermayer to the multivariate setting seems not to be possible. Indeed, Bouchard and Pham comment that ``\emph{it turns out that the one-dimensional argument of Kramkov and Schachermayer does not work directly in our multivariate setting}''. One of the important contributions of the current paper is a novel approach to the variational analysis of the dual problem which allows us to prove, even in a multivariate framework, the existence of a solution to the utility maximization problem under the condition of asymptotic satiability of the value function. The relationship between asymptotic satiability of the value function, and finiteness of the convex dual of the value function is made clear in Proposition \ref{thm:satiable}.

As mentioned above, most optimal investment theorems make the stronger assumption of reasonable asymptotic elasticity on the utility function $U$, or a growth condition on the dual function $U^*$ (the notable exceptions being \cite{KramScha03} and \cite{BoucPham}). We show that these types of assumption are included by our results as follows: In Proposition \ref{thm:oldgrowth} we show that if $U$ is bounded from below on the interior of $\RR_+^d$, multivariate risk averse (see Definition \ref{def:MVRA}) and has reasonable asymptotic elasticity (see Definition \ref{def:RAE}) then $U^*$ satisfies a growth condition (see Definition \ref{def:growthcond}). In Corollary \ref{thm:RAEimplAS}, we show that if $U$ is bounded above, or if $U^*$ satisfies the growth condition then the value function $u$ is asymptotically satiable (which is the hypothesis of this paper). We should point out that multivariate risk aversion is not the same as concavity, and we feel that its importance has been overlooked by the existing literature on multivariate utility maximization. In particular, it appears to be an essential ingredient in the proof of Proposition \ref{thm:oldgrowth}.

There are three standard ways to formulate a dual optimization problem in the utility maximization literature: In terms of martingale measures, their density processes or their Radon-Nikod\'ym derivatives. In all three cases, these control sets are not large enough to contain the dual optimizer, and they need to be enlarged in some way. For example, in \cite{KramScha99} the set of (martingale) density processes is enlarged by including supermartingales as the control processes, and they employ an abstract dual problem which is formulated using random variables which have lost some mass. In \cite{DeelPhamTouz}, the set of Radon-Nikod\'ym derivatives is enlarged, by including random variables which have lost some mass. In this paper, we develop further the approach of \cite{CvitSchaWang,KaraZitk,OwenZitk} by considering the enlarged space of (finitely additive) Euclidean vector measures. The domain of the dual problem is then complete in the relevant topology, and thus contains the dual optimizer. In Example \ref{exa:singularcomponent} we show that this enlargement is necessary by providing an example where the dual minimizer has a non-zero singular component. Our approach makes explicit the ``loss of mass'' experienced by the dual minimizer; in previous work on transaction costs, the dual minimizer corresponds to the countably additive part of our dual minimizer. Our approach is just as powerful as the approach of using a dual control process. Indeed, each finitely additive measure in the domain of our dual problem gives rise to a supermartingale control process (see e.g. \cite[Proposition 2.2]{KaraZitk} for this construction in the univariate case).

There have also been several approaches used in the literature to show the absence of a gap between the optimal primal and dual values. These approaches include using minimax, the Fenchel duality theorem, and the Lagrange duality theorem. In a recent paper \cite{KleiRoge}, Klein and Rogers propose a flexible approach which identifies the dual problem for financial markets with frictions. They guarantee the absence of a duality gap by using minimax, under the assumption of a duality condition which they call (XY). We have chosen to follow instead the approach of \cite{OwenZitk}, using the perfectly suited, and equally powerful Lagrange duality theorem as our weapon of choice (see Proposition \ref{thm:duality} and Theorem \ref{thm:lagrange}). Of course, the minimax, Fenchel duality, and Lagrange duality theorems on non-separable vector spaces are all based upon the the Hahn-Banach theorem in its geometric form, the separating hyperplane theorem.

\vspace\baselineskip

The rest of the paper is structured as follows. In Section 2 we introduce some preliminaries, including the transaction costs framework, and some theory of convex analysis, multivariate utility functions and Euclidean vector measures. In Section \ref{sec:main_results} we prove our main theorems, as described above. In Section \ref{sec:liquidation}, we explain how to relate the formulation of our optimization problem to the liquidation case. In the appendix we present the Lagrange duality theorem, which is used to show that there is no duality gap. The appendix also contains the proofs of some of the auxiliary results from Section 2, which are postponed in order to improve the presentation.

\section{Preliminaries}\label{sec:prelim}

In this section we present all the preliminary concepts and notation which are required for the analysis of the optimization problem. The reader may wish to skip these preliminaries at first, and refer back when necessary. The structure of this section is as follows. In Subsection \ref{sec:bidaskformalism} we recall the transaction costs framework of \cite{CampScha}. In Section \ref{sec:conv_analysis} we introduce some terminology from convex analysis, including dual functionals and their properties. In Subsection \ref{sec:multivariate} we introduce multivariate utility functions and discuss various properties such as asymptotic satiability, reasonable asymptotic elasticity, and multivariate risk aversion. Finally, in Subsection \ref{sec:vector_measures} we collect some facts about Euclidean vector measures, which we use for our formulation of the dual problem.

\subsection{Bid-ask matrix formalism of transaction costs}\label{sec:bidaskformalism}

Let us recall the basic features of the transaction costs model as formalized in \cite{CampScha} (see also \cite{Scha}). In such a model, all agents can trade in $D$ assets according to a random and time varying bid-ask matrix. A $D\times D$ matrix $\Pi=(\pi^{ij})_{1\le i,j\le D}$ is called a \emph{bid-ask matrix} if
(i) $\pi^{ij}>0$ for every $1\le i,j\le D$,
(ii) $\pi^{ii}=1$ for every $1\le i\le D$, and
(iii) $\pi^{ij}\le\pi^{ik}\pi^{kj}$ for every $1\le i,j,k\le D$.
The entry $\pi^{ij}$ denotes the number of units of asset $i$ required to purchase one unit of asset $j$. In other words, $1/\pi^{ji}$ and $\pi^{ij}$ denote, respectively, the bid and ask prices of asset $j$ denominated in asset $i$.

Given a bid-ask matrix $\Pi$, the \emph{solvency cone} $K(\Pi)$ is defined as the convex polyhedral cone in $\RR^D$ generated by the canonical basis vectors $e^i$, $1\le i\le D$ of $\RR^D$, and the vectors $\pi^{ij}e^i-e^j$, $1\le i,j\le D$. The cone $-K(\Pi)$ should be intepreted as those portfolios available at price zero. The (positive) polar cone of $K(\Pi)$ is defined by
\[ K^*(\Pi) = \sets{w\in\RR^D}{\ip vw\ge0,\forall v\in K(\Pi)}. \]

Next, we introduce randomness and time in our model. Let $(\Omega,(\F_t)_{t\in[0,T]},\PP)$ be a filtered probability space satisfying the usual conditions and supporting all processes appearing in this paper. An adapted, c\`adl\`ag process $(\Pi_t)_{t\in[0,T]}$ taking values in the set of bid-ask matrices will be called a \emph{bid-ask process}. A bid-ask process $(\Pi_t)_{t\in[0,T]}$ will now be fixed, and we drop it from the notation by writing $K_\tau$ instead of $K(\Pi_\tau)$ for a stopping time $\tau$.

In accordance with the framework developed in \cite{CampScha} we make the following technical assumption throughout the paper. The assumption is equivalent to disallowing a final trade at time $T$, but it can be relaxed via a slight modification of the model (see \cite[Remark 4.2]{CampScha}). For this reason, we shall not explicitly mention the assumption anywhere.
\begin{assumption}\label{ass:contPi}
$\F_{T-}=\F_T$ and $\Pi_{T-}=\Pi_T$ a.s.
\end{assumption}

\begin{definition}
An adapted, $\RR_+^D\setminus\{0\}$-valued, c\`adl\`ag martingale $Z=(Z_t)_{t\in[0,T]}$ is called a \emph{consistent price process} for the bid-ask process $(\Pi_t)_{t\in[0,T]}$ if $Z_t\in K_t^*$ a.s. for every $t\in[0,T]$. Moreover, $Z$ will be called a \emph{strictly consistent price process} if it satisfies the following additional condition: For every $[0,T]\cup\set{\infty}$-valued stopping time $\tau$, $Z_\tau\in\interior(K_\tau^*)$ a.s. on $\set{\tau<\infty}$, and for every predictable $[0,T]\cup\set{\infty}$-valued stopping time $\sigma$, $Z_{\sigma-}\in\interior(K_{\sigma-}^*)$ a.s. on $\set{\sigma<\infty}$. The set of all (strictly) consistent price processes will be denoted by $\Z$ ($\Z^s$).
\end{definition}

The following assumption, which is used extensively in \cite{CampScha}, will also hold throughout the paper.
\begin{assumption}[SCPS]\label{ass:SCPS}
Existence of a strictly consistent price system: $\Z^s\neq\emptyset$.
\end{assumption}
This assumption is intimately related to the absence of arbitrage (see also \cite{JouiKall,GuasRasoScha,GuasRaso}).

\begin{definition}\label{def:admV}
Suppose that $(\Pi_t)_{t\in[0,T]}$ is a bid-ask process such that Assumption \ref{ass:SCPS} holds true. An $\RR^D$-valued process $V = (V_t)_{t\in[0,T]}$ is called a self-financing portfolio process for the bid-ask process $(\Pi_t)_{t\in[0,T]}$ if it satisfies the following properties:
\begin{itemize}
 \item[(i)] It is predictable and a.e. path has finite variation (not necessarily right-continuous).
 \item[(ii)] For every pair of stopping times $0\le\sigma\le\tau\le T$, we have
 \[ V_\tau-V_\sigma \in -\overline{\operatorname{conv}}\left(\bigcup_{\sigma\le t< \tau}K_t,0\right)\quad\text{a.s.} \]
\end{itemize}
A self-financing portfolio process $V$ is called \emph{admissible} if it satisfies the additional property
\begin{itemize}
 \item[(iii)] There is a constant $a>0$ such that $V_T+a\one\in K_T$ a.s. and $\ip{V_\tau+a\one}{Z_\tau^s}\ge0$ a.s. for all $[0,T]$-valued stopping times $\tau$ and for every strictly consistent price process $Z^s\in\Z^s$. Here, $\one\in\RR^D$ denotes the vector whose entries are all equal to $1$.
\end{itemize}
\end{definition}
Let $\mathcal A^x$ denote the set of all admissible, self-financing portfolio processes with initial endowment $x\in\RR^D$, and let
\[ \A{x}:=\sets{V_T}{V\in\mathcal A^x} \]
be the set of all contingent claims attainable at time $T$ with initial endowment $x$. Note that $\A{x}=x+\A{0}$ for all $x\in\RR^D$.
\begin{remark}
  A few observations about the previous definition of admissible self-financing strategy are in order. We recall that for any portfolio process $V=(V^1 ,\dots,V^D)$, the quantity $V_t^i$ (for $1\leq i\leq D$) represents the number of units of asset $i$ held by the agent at time $t$. The condition of a.s. finite variation in (i) is justified by the fact that, since for each change in the portfolio the agent must pay a proportional transaction cost, the transaction costs would add up to infinity for trajectories with infinite variation. It has been shown in \cite{GuasRasoScha,GuasRaso} that in a one-dimensional setting this property is a consequence of the assumption of No-Free-Lunch. Therefore it is economically meaningful to restrict to portfolio processes with a.e. trajectory of finite variation. 

  Condition (ii) can be translated in these terms: Fixing stopping times $\sigma\le\tau$, the portfolio's change $V_\tau-V_\sigma$ should be a.s. in the closure of the sum of the cones $(-K_t)_{t\in[\sigma,\tau)}$ of contingent claims available (at time $t$) at price zero. This is the analogue of the self-financing condition usually considered in the frictionless case.
  
  For a more detailed discussion of the content of Definition \ref{def:admV}, especially the very delicate admissibility condition (iii) and the reasons why portfolio processes are allowed to have jumps from the right, we refer to \cite{CampScha}.
\end{remark}

For the convenience of the reader we present a reformulation of \cite[Theorem 4.1]{CampScha}, which will be an essential ingredient in the proof of Theorem \ref{thm:optimal}. 
\begin{theorem}[Super-replication]\label{thm:superrep}
Let $x\in\RR^D$ and let $X$ be an $\F_T$-meas\-ur\-able, $\RR_+^D$-valued random variable. Under Assumption \ref{ass:SCPS} we have
\[ X\in\A{x}\qquad\text{if and only if}\qquad\E{\ip X{Z_T^s}}\le\ip x{Z_0^s}\text{ for all }Z^s\in\Z^s. \]
\end{theorem}

\subsection{Convex analysis}\label{sec:conv_analysis}

Let $(\X,\tau)$ be a locally convex topological vector space, and let $\X^*$ denote its dual space. On the first reading of this section, $\X$ should simply be thought of as Euclidean space $\RR^d$, and $\tau$ the associated Euclidean topology. However, from Section \ref{sec:main_results} onwards we will need the full generality of topological vector spaces. Given a set $S\subseteq\X$ we let $\cl(S)$, $\interior(S)$, $\ri(S)$ and $\aff(S)$ denote respectively the closure, interior, relative interior and affine hull of $S$. We shall say that a set $C\subseteq\X$ is a convex cone if $\lambda C+\mu C\subseteq C$ for all $\lambda,\mu\ge0$. Given set $S\subseteq\X$, we denote its polar cone by
\[ S^* := \sets{x^*\in\X^*}{\ip x{x^*}\ge0\;\forall x\in S}. \]
Note that $S^*$ is weak$^*$ closed. A convex cone $C\subseteq\X$ induces a preorder $\presucc_C$ on $\X$: We say that $x,x'\in\X$ satisfy $x'\presucc_Cx$ if and only if $x'-x\in C$.

Let $\U:\X\rightarrow[-\infty,\infty]$ be a concave functional on $\X$, that is, the hypograph
\[ \hypo(\U):=\sets{(x,\mu)}{x\in\X,\;\mu\in\RR,\; \mu\le \U(x)} \]
is convex as a subset of $\X\times\RR$. The \emph{effective domain}, $\dom(\U)$, of $\U$ is the projection of $\hypo(\U)$ onto $\X$, i.e. $\dom(\U):=\sets{x\in\X}{\U(x)>-\infty}$. The functional $\U$ is said to be \emph{proper concave} if its effective domain is nonempty, and it never assumes the value $+\infty$.

The \emph{closure}, $\cl(\U)$, of the functional $\U$ is the unique functional whose hypograph is the closure of $\hypo(\U)$ in $\X\times\RR$. The functional $\U$ is said to be \emph{closed} if $\cl(\U)=\U$.

The functional $\U$ is said to be \emph{upper semi-continuous} if for each $c\in\RR$ the set $\sets{x\in\X}{\U(x)\ge c}$ is closed. Equivalently, $\U$ is upper semi-continuous if $\limsup_\alpha\U(x_\alpha)\le\U(x)$, whenever $(x_\alpha)_{\alpha\in A}\subseteq \X$ is a net tending to some $x\in\X$. It is an elementary result that a concave functional is closed if and only if it is upper semi-continuous (see e.g. \cite[Theorem 2.2.1]{Zali} or \cite[Corollary 2.60]{AlipBord}).

Let $\partial \U(x)$ denote the superdifferential of $\U$ at $x$. That is, $\partial \U(x)$ is the collection of all $x^*\in\X^*$ such that
\[ \U(z)\le \U(x)+\ip{z-x}{x^*}\qquad\forall z\in\X. \]

A functional $\V:\X\rightarrow[-\infty,\infty]$ is said to be \emph{convex} if $-\V$ is concave. The corresponding definitions of the effective domain, proper convexity, the lower semi-continuity, closure and subdifferential for a convex functional are made in the obvious way. 

\begin{definition}[Dual functionals]\label{def:dualfunctionals}
\begin{enumerate}
\item If $\U:\X\rightarrow[-\infty,\infty)$ is proper concave then we define its dual functional $\U^*:\X^*\rightarrow(-\infty,\infty]$ by
    \begin{equation}\label{eqn:conv_conj}
      \U^*(x^*):=\sup_{x\in\X}\set{\U(x)-\ip x{x^*}}.
    \end{equation}
    The dual functional $\U^*$ is a weak$^*$ lower semi-continuous, proper convex functional on $\X^*$. Note that $\U^*=(\cl(\U))^*$ (see e.g. \cite[Theorem 2.3.1]{Zali}).
\item If $\V:\X^*\rightarrow(-\infty,\infty]$ is proper convex then we define the pre-dual functional ${}^*\V:\X\rightarrow[-\infty,\infty)$ by
    \[ {}^*\V(x):=\inf_{x^*\in\X^*}\set{\V(x^*)+\ip x{x^*}}. \]
    Similarly, ${}^*\V$ is a weakly\footnote{A concave functional is weakly upper semi-continuous if and only if it is originally upper semi-continuous.} 
    upper semi-continuous, proper concave functional. By applying \cite[Theorem 2.3.3]{Zali} we see that $({}^*\V)^*=\cl \V$.
\end{enumerate}
\end{definition}

The reader should be aware that the dual functional is not the same object as the conjugate functional commonly used in texts on convex analysis. Nevertheless the only discrepancies are in the sign convention; any property of conjugate functions can, with a little care, be re-expressed as a property of the dual function.

The next lemma will be used several times throughout the paper. Its proof is simple, and is therefore omitted. We say that $\U$ is \emph{increasing} with respect to a preorder $\presucc$ on $\X$, if $\U(x')\ge \U(x)$ for all $x,x'\in\X$ such that $x'\presucc x$.

\begin{lemma}\label{thm:dual_functional}
Let $\U:\X\rightarrow[-\infty,\infty)$ be proper concave. Then $\U^*$ is decreasing with respect to the preorder induced by $(\dom(\U))^*$. Suppose furthermore that $\U$ is increasing with respect to the preorder induced by some cone $C$. Then $\dom(\U^*)\subseteq C^*$.
\end{lemma}

\subsection{Multivariate utility functions}\label{sec:multivariate}

\begin{definition}\label{def:utility}
We shall say that a proper concave function $U:\RR^d\rightarrow[-\infty,\infty)$ is a \emph{(multivariate) utility function} if
\begin{enumerate}
\item $C_U:=\cl(\dom(U))$ is a convex cone such that $\RR_+^d\subseteq C_U\neq\RR^d$; and
\item $U$ is increasing with respect to the preorder induced $C_U$.
\end{enumerate}
We call $C_U$ the \emph{support} (or support cone) of $U$, and say that $U$ is supported on $C_U$. The dual function $U^*$ of a utility function $U:\RR^d\rightarrow\RR$ is defined by \eqref{eqn:conv_conj}, with $\X=\RR^d$.
\end{definition}

We shall focus on three particular utility functions in this paper: The agent's utility function $U$ is assumed to be supported on $\RR_+^d$, the extended utility function $\tilde U$ defined by \eqref{eqn:extension} is therefore supported on $\RR_+^D$, and we shall show in Proposition \ref{thm:value_function} that under Assumption \ref{ass:proper_concave} the value function $u$ defined by \eqref{eqn:primal} is a utility function which is supported on a cone which is strictly larger than $\RR_+^D$.

\begin{example}\label{exa:utilfuncs}
\begin{enumerate}
 \item The canonical univariate utility functions on $\RR_+$ are constant relative risk aversion (CRRA) utility functions. These are defined, for $x\in\RR_+$, by
 \[ U_\gamma(x) = \begin{cases} x^\gamma/\gamma, &\gamma<1, \gamma\neq0, \\ \ln x+1/2, & \gamma=0, \end{cases} \]
 with $U_\gamma(x)=-\infty$ otherwise. The dual functions are $U_\gamma^*=-U_{\gamma^*}$ where $\gamma^*$ is the conjugate of the elasticity $\gamma$ (that is, $1/\gamma+1/\gamma^*=1$, unless $\gamma=0$, in which case $\gamma^*=0$).

 \item The simplest class of utility functions which are supported on $\RR_+^d$, is the class of additive utility functions,
 \[ U(x_1,\dots,x_d) := \sum_{i=1}^d U_i(x_i), \]
 where $U_1,\dots,U_d:\RR\rightarrow[-\infty,\infty)$ are univariate utility functions on $\RR_+$. In this case the dual function also takes the additive form $U^*(x^*)=\sum_{i=1}^d U_i^*(x_i^*)$.

 \item The Cobb-Douglas utility functions form another class of utility functions supported on $\RR_+^d$. Define
 \[ U(x_1,\dots,x_d) := \begin{cases} \prod_{i=1}^d x_i^{\alpha_i}, & x\in\RR_+^d, \\ -\infty, & \text{otherwise}, \end{cases} \]
 where $\alpha_i\ge0$ are such that $\sum_{i=1}^d\alpha_i<1$.
\end{enumerate}
\end{example}

Note that the dual of the extended function $\tilde U:\RR^D\rightarrow\RR$ is given by
 \begin{equation}
   \tilde U^*(x^*)=\begin{cases}U^*(x_1^*,\dots,x_d^*),&x^*\in\RR_+^D\\+\infty,&\text{otherwise.}\end{cases}
 \end{equation}

In the following subsections we investigate a number of conditions which can be imposed on multivariate utility functions.

\subsubsection{Multivariate Inada conditions: Essential smoothness and asymptotic satiability}

In this subsection we investigate analogues of the well known ``Inada conditions'' for the case of a smooth multivariate utility function. The first condition, which we recall from \cite{Rock}, is well known within the field of convex analysis.
\begin{definition}\label{def:ess_smooth}
A proper concave function $U:\RR^d\rightarrow[-\infty,\infty)$ is said to be \emph{essentially smooth} if
\begin{enumerate}
\item $\interior(\dom(U))$ is nonempty;
\item $U$ is differentiable throughout $\interior(\dom(U))$;
\item $\lim_{i\rightarrow\infty}|\nabla U(x_i)|=+\infty$ whenever $x_1,x_2,\dots$ is a sequence in $\interior(\dom(U))$ converging to a boundary point of $\interior(\dom(U))$.
\end{enumerate}
A proper \emph{convex} function $V$ is said to be essentially smooth if $-V$ is essentially smooth.
\end{definition}

The next result is well known, and can be deduced by a standard application of \cite[Theorems 7.4, 12.2, 26.1, 26.3 and Corollary 23.5.1]{Rock}.
\begin{lemma}\label{thm:essconc_esssmooth}
Let $U$ be a proper concave function which is essentially smooth and strictly concave on $\interior(\dom(U))$. Then $U^*$ is strictly convex on $\interior(\dom(U^*))$, and essentially smooth. Moreover, the maps $\nabla U:\interior(\dom(U))\rightarrow\interior(\dom(U^*))$ and $\nabla U^*:\interior(\dom(U^*))\rightarrow-\interior(\dom(U))$ are bijective and $(\nabla U)^{-1}=-\nabla U^*$.
\end{lemma}

The next condition appears to be less well known.

\begin{definition}\label{def:ass_satiable}
We say that a utility function $U$ is \emph{asymptotically satiable} if for all $\epsilon>0$ there exists an $x\in\RR^d$ such that $\partial(\cl(U))(x)\cap[0,\epsilon)^d\neq\emptyset$.
\end{definition}

The proof of the next lemma can be found in the appendix.
\begin{lemma}\label{thm:nec_suf_ass_sat}
A sufficient condition for asymptotic satiability of $U$ is that for all $\epsilon>0$ there exists an $x\in\interior(\dom(U))$ such that $\partial U(x)\cap[0,\epsilon)^d\neq\emptyset$. If $U$ is either upper semi-continuous or essentially smooth then the condition is both necessary and sufficient for asymptotic satiability.
\end{lemma}

Asymptotic satiability means that one can find positions for which the utility function is almost horizontal. The economic intepretation of this condition is even clearer if $U$ is multivariate risk averse (see Subsection \ref{sec:mvRAandAE}). In this case, the marginals of $U$ decrease with increasing wealth, which means that an asymptotically satiable utility function approaches horizontality in the limit as the quantities of assets consumed increase to infinity.

Let us now consider the effect of asymptotic satiability on the dual function. Recall that for a utility function $U$ we define the closed, convex cone $C_U:=\cl(\dom(U))$. Since the dual function $U^*$ of a utility function is convex, it follows that $\cl(\dom(U^*))$ is convex. Furthermore, as an immediate consequence of Lemma \ref{thm:dual_functional}, we have that $\cl(\dom(U^*))\subseteq (C_U)^*\subseteq\RR_+^d$, and $U^*$ is decreasing with respect to $\presucc_{(C_U)^*}$. However it can happen that $\cl(\dom(U^*))$ fails to be a convex cone, in which case it is strictly contained in $(C_U)^*$. In Proposition \ref{thm:satiable} we give a simple condition under which $\cl(\dom(U^*))=(C_U)^*$. Its proof can be found in the appendix.
\begin{proposition}\label{thm:satiable}
Let $U$ be a utility function. The following conditions are equivalent:
\begin{enumerate}
 \item $U$ is asymptotically satiable;
 \item $0\in\cl(\dom(U^*))$;
 \item $\cl(\dom(U^*))=(C_U)^*$; and
 \item $\cl(\dom(U^*))$ is a convex cone.
\end{enumerate}
If $U$ is asymptotically satiable then we define the closed convex cone $C_{U^*}:=\cl(\dom(U^*))$, so that condition 3 can be written more succinctly as $C_{U^*}=(C_U)^*$.
\end{proposition}

One should think of essential smoothness and asymptotic satiability as the multivariate analogues of the univariate Inada conditions $U'(0)=\infty$ and $U'(\infty)=0$ respectively. Indeed, an additive utility function (see part 2 of Example \ref{exa:utilfuncs}) with continuously differentiable components, $U_i$ ($i=1,\dots,d$), is essentially smooth if and only if each component satisfies $U_i'(0)=\infty$, and asymptotically satiable if and only if each component satisfies $U_i'(\infty)=0$. Clearly these conditions reduce to the usual Inada conditions in the univariate case.

The proof of the following corollary of Lemma \ref{thm:essconc_esssmooth} and Proposition \ref{thm:satiable} is straightforward, and is therefore omitted.

\begin{corollary}\label{thm:functionI}
Let $U:\RR^d\rightarrow[-\infty,\infty)$ be a utility function which is supported on $\RR_+^d$, and which satisfies Assumption \ref{ass:basic}. Recall that by definition of the dual function we have
\begin{equation}\label{eqn:also_no_idea}
 U^*(x^*) \ge U(x)-\ip x{x^*}
\end{equation}
for all $x,x^*\in\RR^d$. If $x^*\in\interior(\RR_+^d)$ then we have equality in \eqref{eqn:also_no_idea} if and only if $x=I(x^*):=-\nabla U^*(x^*)$.

Given $D\ge d$, define $\tilde U:\RR^D\rightarrow[-\infty,\infty)$ by \eqref{eqn:extension}. Again, by definition of the dual function we have
\begin{equation}\label{eqn:Ive_run_out_of_ideas_for_labels}
 \tilde U^*(x^*) \ge\tilde U(x)-\ip x{x^*},
\end{equation}
for all $x,x^*\in\RR^D$. Define $\proj:\RR^D\rightarrow\RR^d$ by   
\begin{equation}\label{eqn:def_proj}
 \proj(x_1,\dots,x_d,x_{d+1},\dots,x_D):=(x_1,\dots,x_d),
\end{equation}
and $\tilde I:\interior(\RR_+^d)\times\RR_+^{D-d}\rightarrow\interior(\RR_+^d)\times\RR_+^{D-d}$ by
\begin{equation}\label{def:functionI} \tilde I(x^*):=(-\nabla U^*(\proj(x^*)),\underline 0), \end{equation}
where $\underline 0$ denotes the zero vector in $\RR^{D-d}$. Then, \textup{(i)} if $x^*\in\interior(\RR_+^d)\times\RR_+^{D-d}$ then we have equality in \eqref{eqn:Ive_run_out_of_ideas_for_labels} whenever $x=\tilde I(x^*)$ and \textup{(ii)} if $x^*\in\interior(\RR_+^D)$ then there is equality in \eqref{eqn:Ive_run_out_of_ideas_for_labels} if and only if $x=\tilde I(x^*)$.
\end{corollary}

\subsubsection{Multivariate risk aversion}

In this subsection we present the multivariate analogue of risk aversion. Generalisation of the concept of risk aversion to the multivariate case was first considered in \cite{Rich}. The idea is that a risk-averse investor should prefer a lottery in which they have an even chance of winning $x+z$ or $x+z'$ (with $z,z'$ positive), to a lottery in which they have an even chance of winning $x$ or $x+z+z'$. Put differently, the investor prefers lotteries where the outcomes are less extreme. Some further, mathematically equivalent conditions for multivariate risk aversion can be found in \cite[Theorem 3.12.2]{MuelStoy}.

In one dimension, multivariate risk aversion is equivalent to concavity of the utility function, however in higher dimensions this is no longer the case.
\begin{definition}\label{def:MVRA}
\begin{enumerate}
 \item Let $U$ be a utility function which is supported on $\RR_+^d$. We shall say that $U$ is \emph{multivariate risk averse} if for any $x\in\RR^d$ and any $z,z'\in\RR_+^d$ we have
 \begin{equation}\label{eqn:MVRA}
   U(x)+U(x+z+z')\le U(x+z)+U(x+z');
 \end{equation}
 \item Let $U$ be a utility function which is supported on $\RR_+^d$. We shall say that $U$ has \emph{decreasing marginals} if for any $x\in\dom(U)$, any $x'\in\RR^d$ satisfying $x'\presucc_{\RR_+^d}x$, and any $z\in\RR_+^d$ we have
 \[ U(x+z)-U(x)\ge U(x'+z)-U(x'). \]
\end{enumerate}
\end{definition}

The proof of the following result is simple, and is therefore omitted.
\begin{lemma}\label{thm:mvra}
Let $U$ be a utility function which is supported on $\RR_+^d$. Then $U$ is multivariate risk aversion if and only if it has decreasing marginals. If $U$ is differentiable on $\interior(\RR_+^d)$ and multivariate risk averse then given $x,x'\in\interior(\RR_+^d)$ such that $x'\presucc_{\RR_+^d}x$ we have $\nabla U(x)\presucc_{\RR_+^d}\nabla U(x')$.
\end{lemma}

If $U$ is an additive utility function (see part 2 of Example \ref{exa:utilfuncs}) then the concavity of each component $U_i$ is enough to imply that $U$ is multivariate risk averse. However not all utility functions are multivariate risk averse; the Cobb-Douglas utility functions (see part 3 of Example \ref{exa:utilfuncs}) provide examples of such utility functions. To get a better feel for why, in the general case, multivariate risk aversion is \emph{not} the same as concavity, it helps to consider the Hessian of a (twice differentiable) utility function. The utility function exhibits multivariate risk aversion if at every point the Hessian contains only non-positive entries; in other words, all second order partial derivatives are non-positive. In contrast, the Hessian of a concave function at every point is negative semi-definite.

\subsubsection{Reasonable asymptotic elasticity and the growth condition}\label{sec:mvRAandAE}

We begin by presenting a multivariate analogue of the well known condition of reasonable asymptotic elasticity.
\begin{definition}\label{def:RAE}
Let $U$ be an essentially smooth utility function which is supported on $\RR_+^d$, and bounded from below on $\interior(\RR_+^d)$. We say that $U$ has \emph{reasonable asymptotic elasticity} if
\begin{equation}\label{eqn:oldasymptelast}
\sup_{c\in\RR}\liminf_{\substack{x\in\interior(\RR_+^d)\\|x|\rightarrow\infty}}\frac{U(x)+c}{\ip x{\nabla U(x)}}>1,
\end{equation}
where $|x|:=\max\set{|x_1|,\dots,|x_d|}$.
\end{definition}

As an example, the additive utility function $U(x)=\sum_{i=1}^d U_i(x_i)$, with $U_i(x_i):=x_i^{\gamma_i}/\gamma_i$, $x_i>0$, where $0<\gamma_i<1$ for each $i=1,\dots,d$ (see part 2 of Example \ref{exa:utilfuncs}) has reasonable asymptotic elasticity.

The definition of asymptotic elasticity in the univariate setting is due to \cite{KramScha99}. In the multivariate setting, one can define the asymptotic elasticity of an essentially smooth utility function supported on $\RR_+^d$ by
\begin{equation}\label{eqn:AE}
\operatorname{AE}(U):=\limsup\sets{\ip x{\nabla U(x)}/U(x)}{x\in\interior(\RR_+^d),\;|x|\rightarrow\infty},
\end{equation}
provided the utility function $U$ is strictly positive on $\interior(\RR_+^d)$. In this case, it is trivial that if $\operatorname{AE}(U)<1$ then \eqref{eqn:oldasymptelast} holds. We prefer to formulate the condition of reasonable asymptotic elasticity in terms of the reciprocal of the ratio used in \eqref{eqn:AE}, since the term $\ip x{\nabla U(x)}$ in the denominator of \eqref{eqn:oldasymptelast} is guaranteed to be strictly positive for all $x\in\interior(\RR_+^d)$. Note that the assumption in equation \eqref{eqn:AE}, that $U$ is strictly positive on $\interior(\RR_+^d)$, is relaxed in Definition \ref{def:RAE} to allow $U$ which are bounded below on $\interior(\RR_+^d)$, effectively by adding the constant $c$. Note also that the supremum in \eqref{eqn:oldasymptelast} can be replaced by the limit as $c\rightarrow\infty$.

Unfortunately it is senseless to extend Definition \ref{def:RAE} to the case where $U$ is unbounded below on $\interior(\RR_+^d)$, unless $d=1$. Indeed, by inspection of \eqref{eqn:oldasymptelast}, it is clear that a necessary condition for a utility function to have reasonable asymptotic elasticity is the existence of a sublevel set $\{x\in\interior(\RR_+^d):U(x)\le-c\}$ which is either bounded or empty, a condition which fails whenever $d\ge2$ for additive utility functions which are unbounded from below on $\interior(\RR_+^d)$.

Variations of Definition \ref{def:RAE} have already appeared in the literature for the case where $U(0)=0$ and $U(\infty)=\infty$ (see e.g. \cite{DeelPhamTouz,Kami01,Kami04}). At a first glance, the differences between the definitions of reasonable asymptotic elasticity in these three papers appear to be slight, however more thought reveals that this is in fact a rather delicate issue.

In each of the three papers mentioned, the assumption of reasonable asymptotic elasticity is used in order to prove a growth condition on the dual function $U^*$ (see Definition \ref{def:growthcond}). In turn, the growth condition can be used as an ingredient in the proof of the existence of the optimizer in the primal problem. However, it appears that the definitions of reasonable asymptotic elasticity in \cite{DeelPhamTouz} and \cite{Kami01} are not quite strong enough to imply the growth condition. To compensate for this, Kamizono uses, for instance, an additional assumption (4.22b) which unfortunately excludes all additive utility functions.

Our definition of reasonable asymptotic elasticity is essentially equivalent to the one used in \cite{Kami04}. However, in order to prove the growth condition, we believe the additional assumption of multivariate risk aversion is an essential ingredient (see Proposition \ref{thm:oldgrowth}).

\begin{definition}\label{def:growthcond}
Let $U:\RR^d\rightarrow[-\infty,\infty)$ be a utility function which is supported on $\RR_+^d$, and which is asymptotically satiable. We shall say that the dual function $U^*$ satisfies the \emph{growth condition} if there exists a function $\zeta:(0,1]\rightarrow[0,\infty)$ such that for all $\epsilon\in(0,1]$ and all $x^*\in\interior(\RR_+^d)$
\begin{equation}\label{eqn:growth}
U^*(\epsilon x^*)\le \zeta(\epsilon)(U^*(x^*)^++1).
\end{equation}
\end{definition}

\begin{remark}\label{rem:something}
If $U$ is bounded from above then $U^*$ trivially satisfies the growth condition with $\zeta(\epsilon):=\sup_{x^*\in\RR_+^d}U^*(x^*)=U^*(0)=\sup_{x\in\RR^d}U(x)<\infty$. As an example, if $U(x)=\sum_{i=1}^d U_i(x_i)$ is an additive utility function with $U_i(x_i)=\alpha_ix_i^{\gamma_i}/\gamma_i$, where $\alpha_i>0$ and $\gamma_i<0$ for each $i=1,\dots,d$ (see part 2 of Example \ref{exa:utilfuncs}) then $U^*$ trivially satisfies the growth condition.
\end{remark}

The following two results shed further light on the relationship between the condition of reasonable asymptotic elasticity and the growth condition. Their proofs are provided in the appendix.

\begin{proposition}\label{thm:oldgrowth}
Let $U$ be a utility function which is supported on $\RR_+^d$, and which satisfies Assumption \ref{ass:basic}. If $U$ is bounded from below on $\interior(\RR_+^d)$, multivariate risk averse, and reasonably asymptotically elastic then $U^*$ satisfies the growth condition.
\end{proposition}

\begin{lemma}\label{thm:addUgrowth}
Let $U(x)=\sum_{i=1}^dU_i(x_i)$ be an additive utility function (supported on $\RR_+^d$), which is bounded from below on $\interior(\RR_+^d)$. If each of the components, $U_i$, has reasonable asymptotic elasticity then $U^*$ will satisfy the growth condition.
\end{lemma}

If a utility function is unbounded below on $\interior(\RR_+^d)$ then the previous two results do not apply. It seems therefore that if the utility function is unbounded above and below (on $\interior(\RR_+^d)$) then the growth condition has to be verified on a case-by-case basis. For example, if $U(x_1,x_2):=\ln x_1+\ln x_2+1$ then $U^*$ satisfies the growth condition, while if $U(x_1,x_2)=2x_1^{1/2}-x_2^{-1}$ then $U^*$ fails to satisfy the growth condition.

\subsection{Euclidean vector measures}\label{sec:vector_measures}

A function $m$ from a field $\F$ of subsets of a set $\Omega$ to a Banach space $\X$ is called a finitely additive vector measure, or simply a vector measure if $m(A_1\cup A_2)=m(A_1)+m(A_2)$, whenever $A_1$ and $A_2$ are disjoint members of $\F$. The theory of vector measures was heavily developed in the late 60s and early 70s, and a survey of this theory can be found in \cite{DiesUhl}. In this paper, we will be concerned with the special case where $\X=\RR^D$; we refer to the associated vector measure as a ``Euclidean vector measure'', or simply a ``Euclidean measure''. In this setting, many of the subtleties of the general Banach space theory do not appear. For instance, there is no distinction between the properties of boundedness, boundedness in (total) variation, boundedness in semivariation and strong boundedness. In fact, we can obtain all the results that we need about Euclidean measures by decomposing them into their one-dimensional components. For this reason, we appeal exclusively to results of \cite{RaoRao}, which covers the one-dimensional case very thoroughly.

Let us recall a few definitions from the classical, one-dimensional setting. The \emph{total variation} of a (finitely additive) measure $m:\F\rightarrow\RR$ is the function $|m|:\F\rightarrow[0,\infty]$ defined by
\[ |m|(A):=\sup\sum_{j=1}^n|m(A_j)|, \]
where the supremum is taken over all finite sequences $(A_j)_{j=1}^n$ of disjoint sets in $\F$ with $A_j\subseteq A$. A measure $m$ is said to have \emph{bounded total variation} if $|m|(\Omega)<\infty$. A measure $m$ is said to be \emph{bounded} if $\sup\sets{|m(A)|}{A\in\F}<\infty$. It is straightforward to show that
\[ \sup\sets{|m(A)|}{A\in\F}\le |m|(\Omega)\le2\sup\sets{|m(A)|}{A\in\F}, \]
hence a measure is bounded if and only if it has bounded total variation. A measure $m$ is said to be \emph{purely finitely additive} if $0\le\mu\le|m|$ and $\mu$ is countably additive imply that $\mu=0$. A measure $m$ is said to be \emph{weakly absolutely continuous} with respect to $\PP$ if $m(A)=0$ whenever $A\in\F$ and $\PP(A)=0$.

We turn now to the $D$-dimensional case. A Euclidean measure $m$ can be decomposed into its one-dimensional coordinate measures $m_i:\F\rightarrow\RR$ by defining $m_i(A):=\ip{e^i}{m(A)}$, where $e^i$ is the $i$-th canonical basis vector of $\RR^D$. In this way, $m(A)=(m_1(A),\dots,m_D(A))$ for every $A\in\F$. We shall say that a Euclidean measure $m$ is \emph{bounded}, \emph{purely finitely additive} or \emph{weakly absolutely continuous} with respect to $\PP$ if each of its coordinate measures is bounded, purely finitely additive or weakly absolutely continuous with respect to $\PP$.

Let $\ba(\RR^D)=\ba(\Omega,\F_T,\PP;\RR^D)$ denote the vector space of bounded Euclidean measures $m:\F_T\rightarrow\RR^D$, which are weakly absolutely continuous with respect to $\PP$. Let $\ca(\RR^D)$ the subspace of countably additive members of $\ba(\RR^D)$. Equipped with the norm
\[\|m\|_{\ba(\RR^D)}:=\sum_{i=1}^D|m_i|(\Omega), \]
the spaces $\ba(\RR^D)$ and $\ca(\RR^D)$ are Banach spaces.

Let $\ba(\RR_+^D)$ denote the convex cone of $\RR_+^D$-valued measures within $\ba(\RR^D)$. The next proposition can be easily deduced from its one-dimensional version (see, e.g., \cite[Theorem 10.2.1]{RaoRao}) via a coordinatewise reasoning. Its proof, which also involves a simple application of \cite[Theorems 2.2.1(5), 2.2.2, 10.2.2 and Corollary 10.1.4]{RaoRao}, is therefore omitted.

\begin{proposition}\label{thm:brooks}
Given any $m\in\ba(\RR^D)$ there exists a unique Yosida-Hewitt decomposition $m=m^c+m^p$ where $m^c\in\ca(\RR^D)$ and $m^p$ is purely finitely additive. If $m\in\ba(\RR_+^D)$ then $m^c,m^p\in\ba(\RR_+^D)$.
\end{proposition}

We shall see now that elements of $\ba(\RR^D)$ play a natural role as linear functionals on spaces of (essentially) bounded $\RR^D$-valued random variables. First, some more notation: Let $L^0(\RR^D)=L^0(\Omega,\F_T,\PP;\RR^D)$ denote the space of $\RR^D$-valued random variables (identified under the equivalence relation of a.s. equality). Given $X\in L^0(\RR^D)$ we define the coordinate random variables $X_i\in L^0(\RR)$ for $i=1,\dots,D$ by $X_i:=\ip X{e^i}$, so that $X=(X_1,\dots,X_D)$. Let $L^1(\RR^D)$ denote the subspace of $L^0(\RR^D)$ consisting of those random variables $X$ for which $\|X\|_1:=\E{\sum_i|X_i|}<\infty$. Let $L^\infty(\RR^D)$ denote the subspace of $L^0(\RR^D)$ consisting of those random variables $X$ for which $\|X\|_\infty:=\esssup\big\{\max_i|X_i|\big\}<\infty$. Finally, let $L^\infty(\RR^D)^*$ denote the dual space of $(L^\infty(\RR^D),\|.\|_{\infty})$.

We now define the map $\Psi:\ba(\RR^D)\rightarrow L^\infty(\RR^D)^*$ by
\begin{equation}\label{eqn:multidimensional1}
\big(\Psi(m)\big)(X):=\int_\Omega\ip X{\d m}:=\sum_{i=1}^D\int_\Omega X_i\d m_i,
\end{equation}
where $(m_1,\dots,m_D)$ is the coordinate-wise representation of $m$. For details concerning the construction of the one-dimensional integrals in \eqref{eqn:multidimensional1}, see \cite[Chapter 4]{RaoRao}, where the integral is referred to as the D-integral. We also define the map $\Phi:\ca(\RR^D)\rightarrow L^1(\RR^D)$ by
$\Phi(m):=\left(\radnik[m_1],\dots,\radnik[m_D]\right)$,
where $\radnik[m_i]$ is the Radon-Nikod\'ym derivative of the $i$-th coordinate measure. Finally, we define the isometric embedding $\iota:L^1(\RR^D)\rightarrow L^\infty(\RR^D)^*$ by $\big(\iota(Y)\big)(X):=\E{\ip XY}$. The next proposition can be easily deduced from its one-dimensional version (see, e.g., \cite[Theorem 4.7.10]{RaoRao}) via a coordinatewise reasoning. Its proof is therefore omitted.

\begin{proposition}\label{thm:isometric_isomorphisms}
The maps $\Psi$ and $\Phi$ are isometric isomorphisms. Furthermore, $\iota\circ\Phi=\Psi|_{\ca(\RR^D)}$.
\end{proposition}

\begin{corollary}\label{thm:compact_unit}
$(\ba(\RR^D),\|.\|_{\ba(\RR^D)})$ has a $\sigma(\ba(\RR^D),L^\infty(\RR^D))$-compact unit ball.
\end{corollary}

For the remainder of the paper, we shall overload our notation as follows: Given $m\in\ba(\RR^D)$ and $X\in L^\infty(\RR^D)$, we write $m(X)$ as an abbreviation of $\big(\Psi(m)\big)(X)$, and we define $\radnik[m]:=\left(\radnik[m_1],\dots,\radnik[m_D]\right)=\Phi(m)$.

Given $x\in\RR^D$ and $A\in\F_T$ it follows from equation \eqref{eqn:multidimensional1} that $m(x\chi_A^{})=\ip x{m(A)}$, where $\chi_A^{}$ denotes the indicator random variable of $A$. In the special case where $A=\Omega$, we have $m(x)=\ip x{m(\Omega)}$.

Let $L^0(\RR_+^D)$ and $L^\infty(\RR_+^D)$ denote respectively the convex cones of random variables in $L^0(\RR^D)$ and $L^\infty(\RR^D)$ which are $\RR_+^D$-valued a.s. Note that if $m\in\ba(\RR_+^D)$ and $X\in L^\infty(\RR_+^D)$ then $m(X)\ge0$ (see \cite[Theorem 4.4.13]{RaoRao}). This observation allows us to extend the definition of $m(X)$ to cover the case where $m\in\ba(\RR_+^D)$ and $X\in L^0(\RR_+^D)$ by setting
\begin{equation}\label{eqn:goiosi}
m(X):=\sup_{n\in\NN}m\left(X\wedge_{\RR_+^D}(n\one)\right),
\end{equation}
where $\one\in\RR^D$ denotes the vector whose entries are all equal to $1$, and $(x_1,\dots,x_D)\wedge_{\RR_+^D}(y_1,\dots,y_D):=(x_1\wedge y_1,\dots,x_D\wedge y_D)$. It is trivial that \eqref{eqn:goiosi} is consistent with the definition of $m(X)$ for $X\in L^\infty(\RR^D)$. Furthermore, the supremum in \eqref{eqn:goiosi} can be replaced by a limit, because the sequence of numbers is increasing. It follows that given $m_1,m_2\in\ba(\RR_+^D)$, $\lambda_1,\lambda_2,\mu_1,\mu_2\ge0$ and $X_1,X_2\in L^0(\RR_+^D)$, we have
\begin{align*}
  &(\lambda_1m_1+\lambda_2m_2)(\mu_1X_1+\mu_2X_2) \\
  &\qquad\qquad=\lambda_1\mu_1m_1(X_1)+\lambda_1\mu_2m_1(X_2)+\lambda_2\mu_1m_2(X_1)+\lambda_2\mu_2m_2(X_2).
\end{align*}
Note that the final statement of Proposition \ref{thm:isometric_isomorphisms} means that given $m\in\ca(\RR^D)$ and $X\in L^\infty(\RR^D)$ we have $m(X)=\E{\ip X{\radnik[m]}}$. It is easy to show that this property is also true under the extended definition \eqref{eqn:goiosi}.

\section{Main results}\label{sec:main_results}

Throughout this section $U$ denotes a utility function which is supported on $\RR_+^d$. The extension, $\tilde U$, of $U$ to a utility function supported on $\RR_+^D$ is defined by \eqref{eqn:extension}. The value function $u$ is defined by \eqref{eqn:primal}. We shall indicate explicitly where assumptions on the investor's preferences (i.e. Assumptions \ref{ass:basic}, \ref{ass:proper_concave} and \ref{ass:additional}) are used. 

Regarding our model of the economy, Assumptions \ref{ass:contPi} and \ref{ass:SCPS} will be taken as standing assumptions throughout this section. As noted in Subsection \ref{sec:bidaskformalism}, Assumption \ref{ass:contPi} is a technical assumption which can be relaxed, so we shall not mention this assumption anywhere. To avoid mentioning Assumption \ref{ass:SCPS} in the statement of every result, we shall only indicate in the proofs where the assumption is used. As an exception however, we do mention Assumption \ref{ass:SCPS} explicitly in the statement of our main result, Theorem \ref{thm:optimal}.

The following result shows that if $u$ is finite anywhere in the interior of its effective domain, then it is a utility function, and we give a characterization of the closure of the effective domain of $u$.
\begin{proposition}\label{thm:value_function}
Under Assumption \ref{ass:proper_concave} the value function $u$ is a utility function with support cone $C_u:=\cl(\dom(u))=-\{x\in\RR^D:x\in\A0\}$.
\end{proposition}

\begin{proof}
Note first that $u$ is both concave and increasing with respect to $\RR_+^D$, because $\A0$ is convex and $\tilde U$ is both concave and increasing with respect to $\RR_+^D$. We break the proof into the following four steps. We show that (i) $u(x)<\infty$ for all $x\in\RR^D$, (ii) $C_u=-\{x\in\RR^D:x\in\A0\}$, (iii) $C_u\neq\RR^D$ and (iv) $u$ is increasing with respect to $\presucc_{C_u}$.

(i) Suppose, for a contradiction, that there exists some $\tilde x\in\RR^D$ such that $u(\tilde x)=\infty$. By Assumption \ref{ass:proper_concave} there exists an $x\in\interior(\dom(u))$ such that $u(x)<\infty$. Let $a>0$ be large enough so that $x_1:=x+a\one\presucc_{\RR_+^D}\tilde x$. Since $u$ is increasing with respect to $\RR_+^D$, this implies that $u(x_1)\ge u(\tilde x)=\infty$.

Since $x\in\interior(\dom(u))$, there exists an $\epsilon>0$ such that $x_0:=x-\epsilon\one\in\interior(\dom(u))$. We claim that $u(x_0)\in\RR$. Indeed, since $x_0\in\dom(u)$ we have that $u(x_0)>-\infty$, and since $u$ is increasing with respect to $\RR_+^D$, we have $u(x_0)\le u(x)<\infty$.

Since $u(x_0)\in\RR$, we may find an $X_0\in\A{x_0}$ such that $\E{\smash{\tilde U(X_0)}}\in\RR$. Since $u(x_1)=\infty$, given any $R\in\RR$ we may find an $X_1\in\A{x_1}$ such that $\E{\smash{\tilde U(X_1)}}\ge R$. Define now  $\lambda:=\epsilon/(a+\epsilon)\in(0,1)$ and $X:=(1-\lambda)X_0+\lambda X_1\in\A{(1-\lambda)x_0+\lambda x_1}=\A x$. Since $\tilde U$ is concave,
\begin{align*}
 u(x) &\ge \E{\tilde U(X)}=\E{\tilde U((1-\lambda)X_0+\lambda X_1)} \\
 &\ge (1-\lambda)\E{\tilde U(X_0)}+\lambda\E{\tilde U(X_1)}\ge(1-\lambda)\E{\tilde U(X_0)}+\lambda R.
\end{align*}
Since $R$ can be chosen arbitrarily large, this implies that $u(x)=\infty$, which is the required contradiction.

(ii) The set $C:=\{x\in\RR^D:x\in\A0\}$ is a convex cone in $\RR^D$. It follows immediately from \cite[Theorem 3.5]{CampScha} (which requires Assumption \ref{ass:SCPS}) that $C$ is closed in $\RR^D$. Take $x\in\interior(C)$. There exists $\epsilon>0$ such that $x+\epsilon\one\in C$ and hence $\epsilon\one\in\A{-x}$. Now $u(-x)\ge\E{\smash{\tilde U}(\epsilon\one)}=\tilde U(\epsilon\one)>-\infty$, so $-x\in\dom(u)$.

Suppose now that $x\in\dom(u)$. Then $\A x\cap L^0(\RR_+^D)\neq\emptyset$, otherwise this would contradict $u(x)>-\infty$. Pick any $X\in\A x\cap L^0(\RR_+^D)$. Since we may write $0=X-X\in\A x-L^0(\RR_+^D)$ it follows that $0\in\A x$, and hence $x\in -C$.

Since $C$ is closed and $-\interior(C)\subseteq\dom(u)\subseteq-C$, we have $C_u=\cl(\dom(u))=-C$.

(iii) By part (ii), it suffices to show that $\sets{x\in\RR_+^D}{x\in\A0}=\set0$. To show this, suppose that $x\in\RR_+^D$ satisfies $x\in\A0$. Then there exists an admissible portfolio $V$ such that $V_0=0$ and $V_T=x$. Let $Z^s$ be a strictly consistent price process (such a process exists by Assumption \ref{ass:SCPS}). By \cite[Lemma 2.8]{CampScha}, $\bigip{V_t}{Z_t^s}$ is a super-martingale. Hence $0\le\E{\ip x{Z_T^s}}=\E{\ip{V_T}{Z_T^s}}\le\E{\ip{V_0}{Z_0^s}}=0$, and so $x=0$.

(iv) Take $x\in\RR^D$ and $w\in C_u$. Since, by step (i), $u(x)<\infty$, given any $\epsilon>0$ there exists an $X\in\A{x}$ such that $\E{\tilde U(X)}\ge u(x)-\epsilon$. By step (ii), $0\in\A w$, so $X\in\A{x+w}$. Thus
\[ u(x+w)\ge\E{\tilde U(X)}\ge u(x)-\epsilon. \]
Since $\epsilon>0$ is arbitrary, this implies that $u(x+w)\ge u(x)$.\qed
\end{proof}

The following simple example shows that the value function $u$ can fail to be strictly concave on $\interior(\dom(u))$, and may even fail to be differentiable on $\interior(\dom(u))$.
\begin{example}\label{exa:nonstrictC1}
  Consider the case with $D=2$, where the bid-ask process is given by the deterministic, constant matrix
  \[ \Pi_t:=\begin{pmatrix}1&2\\2&1\end{pmatrix}. \]
  In this case, the solvency cones $K_t\equiv K$ are constant, and generated by the vectors $2e^1-e^2$ and $2e^2-e^1$.
  \begin{enumerate}
  \item With $d=2$, we define $U(x_1,x_2):=\ln x_1+\ln x_2$. It is easy to verify that the value function in this case is
  \[ u(x):=\max_{c\in K} U(x-c)=\begin{cases} 2\ln(2x_1+x_2)-3\ln 2,& x_2>2|x_1|, \\ \ln x_1+\ln x_2, & x_1>0, x_1\le2x_2\le4x_1,\\ 2\ln(x_1+2x_2)-3\ln2,& x_1>0,-x_1<2x_2<x_1, \\  -\infty, & \text{otherwise},\end{cases} \]
  which fails to be strictly concave on $\interior(K)$, but which is differentiable throughout $\interior(K)$.

  \item With $d=1$, we define $U(x):=\ln x$, and define $\tilde U:\RR^2\rightarrow[-\infty,\infty)$ by \eqref{eqn:extension}.
  It is easy to verify that the value function in this case is
  \[ u(x):=\max_{c\in K}\tilde U(x-c)=\begin{cases}\ln(x_1+\frac12x_2),&x_2>\max\set{0,-2x_1}, \\ \ln(x_1+2x_2),& x_1>0,-x_1<2x_2\le0, \\ -\infty & \text{otherwise}, \end{cases} \]
  which is fails to be strictly concave on $\interior(K)$, and fails to be differentiable anywhere along the half line $x_1>0,x_2=0$.
\end{enumerate}
\end{example}

Given any initial portfolio $x\in\RR^D$, we define the proper concave functional $\UU_x:L^\infty(\RR^D)\rightarrow[-\infty,\infty)$ by
\begin{equation}\label{def:UU}
\UU_x(X)=\E{\tilde U(x+X)}.
\end{equation}
Since $\tilde U$ is a utility function which is supported on $\RR_+^D$, $\UU_x$ is increasing with respect to the preorder induced by the convex cone $L^\infty(\RR_+^D)$ and $\dom(\UU_0)\subseteq L^\infty(\RR_+^D)$. Let $\UU_x^*:\ba(\RR^D)\rightarrow(-\infty,\infty]$ denote the dual functional defined by \eqref{eqn:conv_conj}. The dual functional is used directly in our formulation of a dual optimization problem (see equation \eqref{eqn:dual} and Proposition \ref{thm:duality}). The following lemma provides a representation of $\UU_x^*$ in terms of the dual function $\tilde U^*$.

\begin{lemma}\label{thm:repVV}
For any $x\in\RR^D$ we have
\[ \UU_x^*(m)=\begin{cases} \displaystyle\E{\tilde U^*\left(\radnik[m^c]\right)}+m(x) & m\in \ba(\RR^D_+) \\ \infty & \text{otherwise.}
\end{cases} \]
\end{lemma}

\begin{proof}
It suffices to consider the case $x=0$ because, setting $\tilde X:=X+x$,
\begin{align*}
  \UU_x^*(m) &= \sup_{X\in L^\infty(\RR^D)}\{\UU_x(X)-m(X)\}=\sup_{\tilde X\in L^\infty(\RR^D)}\{\UU_0(\tilde X)-m(\tilde X)+m(x)\} \\
  &=\UU_0^*(m)+m(x).
\end{align*}

Since $\UU_0$ is increasing with respect to the preorder induced by $L^\infty(\RR_+^D)$, an application of Lemma \ref{thm:dual_functional} gives that $\dom(\UU_0^*)\subseteq L^\infty(\RR_+^D)^*=\ba(\RR_+^D)$. Take $m\in\ca(\RR^D)$. Then by Proposition \ref{thm:isometric_isomorphisms},
\begin{align*}
  \UU_0^*(m) &= \sup_{X\in L^\infty(\RR^D)}\set{\UU_0(X)-m(X)} = \sup_{X\in L^\infty(\RR^D)}\set{\E{\tilde U(X)-\ip X{\radnik[m]}}} \\
  &\le\E{\tilde U^*\left(\radnik[m]\right)}.
\end{align*}
We show that the last inequality also holds in reverse. For each $n\ge1$ define $\tilde U_n^*:\RR^D\rightarrow\RR$ and $I_n:\RR^D\twoheadrightarrow[0,n]^D$ by
\begin{align*}
\tilde U_n^*(x^*) &:= \max\sets{\tilde U(x)-\ip x{x^*}}{x\prepreq_{\RR_+^D}n\one}, \\
I_n(x^*)   &:= \operatorname{argmax}\sets{\tilde U(x)-\ip x{x^*}}{x\prepreq_{\RR_+^D}n\one}.
\end{align*}
For fixed $x^*\in\RR^D$, the sequence $(\tilde U_n^*(x^*))_{n\ge1}$ is monotone increasing to $\tilde U^*(x^*)$, and the random variable $\tilde U_1^*(\radnik[m])$ is integrable. Using the definition of $\UU_0^*$ and the monotone convergence theorem we have
\begin{align*}
\UU_0^*(m) &\ge \sup_n\E{\tilde U\left(I_n\bigg(\radnik[m]\bigg)\right)-\ip{I_n\left(\radnik[m]\right)}{\radnik[m]}} \\
&= \sup_n\E{\tilde U^*_n\left(\radnik[m]\right)} = \E{\tilde U^*\left(\radnik[m]\right)}.
\end{align*}
To finish the proof, it suffices to show that for $m\in \ba(\RR_+^D)$ we have $\UU_0^*(m)=\UU_0^*(m^c)$. An application of Lemma \ref{thm:dual_functional} shows that $\UU_0^*$ is decreasing with respect to the preorder induced by $\ba(\RR_+^D)$. By Proposition \ref{thm:brooks}, $m^p\in\ba(\RR_+^D)$, thus $m\presucc_{\ba(\RR_+^D)}m^c$, and hence $\UU_0^*(m)\le\UU_0^*(m^c)$.

To prove this inequality in the other direction, take any $u\in\RR$ such that $u<\UU_0^*(m^c)$, and any $\epsilon>0$. There exists an $X\in L^\infty(\RR_+^D)$ such that $\UU_0(X)-m^c(X)\ge u$.
An application of \cite[Theorem 10.3.2]{RaoRao} and the monotone convergence theorem gives the existence of an $A\in\F_T$ such that $m^p(\Omega\setminus A)=0$ and $\E{(\tilde U(X)-\tilde U(\epsilon\one))\chi_A^{}}<\epsilon$. An application of \cite[Theorem 4.4.13(ix)]{RaoRao} shows that $m^p(X\chi_{\Omega\setminus A}^{})=0$. Define $\tilde X=X\chi_{\Omega\setminus A}^{}+\epsilon\one\chi_A^{}$. Then
\begin{align*}
& \UU_0(X)-m^c(X)-\UU_0(\tilde X)+m(\tilde X) \\
& \qquad= \E{(\tilde U(X)-\tilde U(\epsilon\one))\chi_A^{}}+m^p(X\chi_{\Omega\setminus A}^{})-m^c(X\chi_A^{})+\epsilon m(\one\chi_A^{})\\
& \qquad\le \epsilon+0+0+\epsilon m(\one)
\end{align*}
Thus
\[ \UU_0^*(m)\ge\UU_0(\tilde X)-m(\tilde X)\ge\UU_0(X)-m^c(X)-\epsilon-\epsilon m(\one)\ge u-\epsilon(1+m(\one)). \]
Since $u<\UU_0^*(m^c)$ and $\epsilon>0$ are arbitrary we have $\UU_0^*(m)\ge\UU_0^*(m^c)$.\qed
\end{proof}

\begin{remark}
Measures in $\dom(\UU_0^*)$ are commonly said to have \emph{finite generalized entropy}. Due to the above characterisation of $\UU_x^*$, it's clear that $\dom(\UU_x^*)=\dom(\UU_0^*)$ for any $x\in\RR^D$.
\end{remark}

Define $\C:=\A0\cap L^\infty(\RR^D)$. The dual cone to $\C$ is defined by
\[ \D:=(-\C)^*=\{m\in \ba(\RR^D):m(X)\le0\text{ for all }X\in\C\}. \]
Note that since $-L^\infty(\RR_+^D)\subseteq\C$, we have $\D\subseteq\ba(\RR_+^D)$.

Given any $x\in\RR^D$ it follows from the definitions of $\D$ and $\UU_x^*$ that
\begin{align}\label{eqn:fenchel}
\sup_{X\in\C}\UU_x(X) &\le \sup_{X\in L^\infty(\RR^D)}\inf_{m\in\D}L_x(X,m) \notag\\
&\le \inf_{m\in\D}\sup_{X\in L^\infty(\RR^D)}L_x(X,m)= \inf_{m\in\D}\UU_x^*(m),
\end{align}
where $L_x(X,m):=\UU_x(X)-m(X)$ is a Lagrangian. Inequality \eqref{eqn:fenchel} is known as Fenchel's inequality, and it identifies
\begin{equation}\label{eqn:dual}
\inf\sets{\UU_x^*(m)}{m\in\D}
\end{equation}
as a potential dual optimization problem.

In our next result, we show that there is no duality gap in \eqref{eqn:fenchel} provided the initial portfolio $x$ does not lie on the boundary of $\dom(u)$. We also show that the dual problem has a solution whenever $x$ lies in the interior of $\dom(u)$.

\begin{proposition}[Duality]\label{thm:duality}
Suppose that Assumption \ref{ass:proper_concave} holds.
\begin{enumerate}
\item For any $x\in\RR^D$ we have
    \begin{equation}\label{eqn:ineq}
      \sup_{X\in\C}\UU_x(X)\le u(x)\le\inf_{m\in\D}\UU_x^*(m).\phantom{=-\infty.}
    \end{equation}
\item If $x\in\interior(\dom(u))=\interior(C_u)$ then
    \[ \sup_{X\in\C}\UU_x(X)=u(x)=\min_{m\in\D}\UU_x^*(m)\in\RR.\;\;\,\;\; \]
\item If $x\not\in\cl(\dom(u))=C_u$ then
    \[ \sup_{X\in\C}\UU_x(X)=u(x)=\inf_{m\in\D}\UU_x^*(m)=-\infty. \]
\end{enumerate}
\end{proposition}

\begin{proof}
\begin{enumerate}
\item The left-hand inequality in \eqref{eqn:ineq} follows trivially from the definitions of $\UU_x$, $\C$ and $u$. To prove the right-hand inequality we need to show that $\E{\tilde U(X)}\le\UU_x^*(m)$ for all $X\in\A{x}$ and $m\in\D$. We may assume without loss of generality that $X\in L^0(\RR_+^D)$, otherwise there is nothing to prove. In this case, for each $n\in\NN$ we have $X\wedge_{\RR_+^D}(n\one)-x\in\C$, and hence
\begin{equation}\label{eqn:flyhjf}
m(X) = \sup_{n\in\NN}m\left(X\wedge_{\RR_+^D}(n\one)\right)
       =m(x)+\sup_{n\in\NN}m\left(X\wedge_{\RR_+^D}(n\one)-x\right)
       \le m(x).
\end{equation}
Furthermore, since $m\in\ba(\RR_+^D)$, it follows from Propositions \ref{thm:brooks} and \ref{thm:isometric_isomorphisms} that
\begin{equation}\label{eqn:adghh}
m(X)=m^c(X)+m^p(X)\ge\E{\ip X{\radnik[m^c]}}+0.
\end{equation}
Using the definition of $\tilde U^*$, combined with equations \eqref{eqn:adghh}, \eqref{eqn:flyhjf} and Lemma \ref{thm:repVV} gives
\begin{align}\label{eqn:fenchel3}
\E{\tilde U(X)}&\le\E{\tilde U^*\left(\radnik[m^c]\right)+\ip X{\radnik[m^c]}} \notag\\
&\le\E{\tilde U^*\left(\radnik[m^c]\right)}+m(x)=\UU_x^*(m).
\end{align}

\item Suppose that $x\in\interior(C_u)$. In order to apply the Lagrange Duality Theorem we set $\X=L^\infty(\RR^D)$ and define the concave functional $\U:\X\rightarrow[-\infty,\infty)$ by $\U=\UU_x$. We must first verify that the hypotheses of part 1 of Theorem \ref{thm:lagrange} hold. Since $x\in\interior(C_u)$, there exists an $\epsilon>0$ such that $x-2\epsilon\one\in C_u$. The deterministic random variable $p:=-\epsilon\one$ lies in the interior of $-L^\infty(\RR^D_+)$ and hence in the interior of $\C$. By Proposition \ref{thm:value_function}, we see that $z:=2\epsilon\one-x\in\A0\cap L^\infty(\RR^D)=\C$. Hence $\U(p+z)=\UU_x(\epsilon\one-x)=\tilde U(\epsilon\one)>-\infty$. Since $x\in\interior(C_u)\subseteq\dom(u)$, part 1 of this proposition gives
\[ \sup_{X\in\C}\U(X)=\sup_{X\in\C}\UU_x(X)\le u(x)<\infty. \]
This verifies the hypotheses of part 1 of Theorem \ref{thm:lagrange}, hence we may assert that
\begin{equation*}
\sup_{X\in \C}\UU_x(X)=\min_{m\in\D}\UU_x^*(m)\in\RR.
\end{equation*}

\item Suppose that $x\not\in C_u$. We set $\X=L^\infty(\RR^D)$ and define the concave functional $\U:\X\rightarrow[-\infty,\infty)$ by $\U=\UU_x$. We must verify that the hypotheses of part 2. of Theorem \ref{thm:lagrange} hold. Since $C_u$ is closed and $x\not\in C_u$, there exists an $\epsilon>0$ such that $x+\epsilon\one\not\in C_u$. The deterministic random variable $p:=-\epsilon\one$ lies in the interior of $\C$. By definition of $C_u$ we have $x-p\not\in\dom(u)$. Using part 1 of this proposition, we see that for any $X\in\C$, $\U(X-p)=\UU_x(X-p)=\UU_{x-p}(X)\le u(x-p)=-\infty$.

By taking any $x'$ in the nonempty set $\interior(C_u)$ and applying part 2. of this proposition, we find the existence of a $\hat m\in\D$ such that $u(x')=\UU_{x'}^*(\hat m)$. Thus by Lemma \ref{thm:repVV}, $\U^*(\hat m)=\UU_x^*(\hat m)=\UU_{x'}^*(\hat m)+\hat m(x-x')=u(x')+\hat m(x-x')<\infty$. This verifies the hypotheses of part 2. of Theorem \ref{thm:lagrange}, and hence we may assert that
\[ \sup_{X\in\C}\UU_x(X)=\inf_{m\in\D}\UU_x^*(m)=-\infty. \tag*{\qed} \]
\end{enumerate}
\end{proof}

The following result will be used in the proofs of Corollary \ref{thm:RAEimplAS} and Proposition \ref{thm:variational}.

\begin{proposition}\label{thm:corol}
Suppose that Assumption \ref{ass:proper_concave} holds. For all $x^*\in\RR^D$ we have
\[ u^*(x^*) = \min\sets{\UU_0^*(m)}{m\in\D\text{ and }m(\Omega)=x^*}, \]
in the sense that the minimum is attained whenever $u^*(x^*)<\infty$.
\end{proposition}

\begin{proof}
Let $v:\RR^D\rightarrow(-\infty,\infty]$ be defined by $v(x^*):=\inf\sets{\UU_0^*(m)}{m\in \D\cap S(x^*)}$, where $S(x^*):=\sets{m\in\ba(\RR_+^D)}{m(\Omega)=x^*}$ and we use the convention that $v(x^*)=\infty$ whenever $\D\cap S(x^*)=\emptyset$.

We begin by showing that the infimum in the definition of $v(x^*)$ is attained whenever $v(x^*)<\infty$. We may assume without loss of generality that $x^*\in\RR_+^D$, otherwise $S(x^*)=\emptyset$. It is straightforward to verify that $S(x^*)$ is a weak$^*$ closed subset of the ball in $\ba(\RR^D)$ of radius $|x^*|_1:=\sum_{i=1}^D|x_i^*|$, and therefore, by Corollary \ref{thm:compact_unit}, $S(x^*)$ is weak$^*$ compact. Since the polar cone $\D$ is weak$^*$ closed this implies that $\D\cap S(x^*)$ is weak$^*$ compact. Since the dual functional $\UU_0^*$ is weak$^*$ lower semi-continuous, the infimum of $\UU_0^*$ over $\D\cap S(x^*)$ is attained whenever $v(x^*)<\infty$.

We claim that $v$ is proper convex. Convexity follows easily from convexity of $\UU_0^*$ and $\D$. That $v$ is proper convex follows from Assumption \ref{ass:proper_concave}, part 2 of Proposition \ref{thm:duality}, Lemma \ref{thm:repVV}, the fact that $\UU_0^*$ is proper convex, and that the minimum in the definition of $v(x^*)$ is attained whenever $v(x^*)<\infty$.

We claim that $v$ is lower semi-continuous. Indeed, suppose that $(x_n^*)_{n\in\NN}\subseteq\RR^D$ is such that $x_n^*\rightarrow x^*$. We may assume without loss of generality that $\liminf_{n\rightarrow\infty}v(x_n^*)<\infty$ otherwise there is nothing to show. There exists a subsequence $(x_{n_k})_{k\in\NN}$ such that $v(x_{n_k}^*)<\infty$ for all $k$, and $\lim_{k\rightarrow\infty}v(x_{n_k}^*)=\liminf_{n\rightarrow\infty} v(x_n^*)$. Let $(\hat m_k)_{k\in\NN}\subseteq\D$ be such that $\hat m_k\in S(x_{n_k}^*)$ and $\UU_0^*(\hat m_k)=v(x_{n_k}^*)$ for each $k$. The sequence $(\hat m_k)_{k\in\NN}$ is bounded in $\ba(\RR^D)$ because for each $k\in\NN$, $\|\hat m_k\|_{\ba(\RR^D)}=|\hat m_k(\Omega)|_1=|x_{n_k}^*|_1\le\sup_{n\in\NN}|x_n^*|_1<\infty$. By Corollary \ref{thm:compact_unit} the sequence $(\hat m_k)_{k\in\NN}$ has a cluster point. There exists, therefore, a directed set $A$, an $\hat m\in\D$ and a subnet $(\hat m_\alpha)_{\alpha\in A}$ of $(\hat m_k)_{k\in\NN}$ which weak$^*$ converges to $\hat m$. Define $x_\alpha^*:=\hat m_\alpha(\Omega)$. The net $(x_\alpha^*)_{\alpha\in A}$ converges to $x^*$. Note that $\hat m\in S(x^*)$ because for each $i=1,\dots,d$ we have $\ip{e^i}{\hat m(\Omega)}=\hat m(e^i)=\lim_\alpha\hat m_\alpha(e^i)=\lim_\alpha\ip{e^i}{\hat m_\alpha(\Omega)}=\lim_\alpha\ip{e^i}{x_\alpha^*}=\ip{e^i}{x^*}$. Since $\hat m\in\D\cap S(x^*)$ and $\UU_0^*$ is weak$^*$ lower semi-continuous, we have $v(x^*)\le\UU_0^*(\hat m)\le\liminf_\alpha\UU_0^*(\hat m_\alpha)=\liminf_\alpha v(x_\alpha^*)=\lim_\alpha v(x_\alpha^*)=\lim_{k\rightarrow\infty} v(x_{n_k}^*)=\liminf_{n\rightarrow\infty}v(x_n^*)$.

By part 2 of Proposition \ref{thm:duality}, and Lemma \ref{thm:repVV} we have, for any $x\in\interior(C_u)$,
\begin{align*}
 u(x) &= \min_{m\in\D}\UU_x^*(m)=\min_{m\in\D}\{\UU_0^*(m)+m(x)\} \\
 &= \min_{x^*\in\RR^D}\min_{\substack{m\in\D\\m(\Omega)=x^*}}\{\UU_0^*(m)+\ip x{x^*}\} \\
 &= \min_{x^*\in\RR^D}\{v(x^*)+\ip x{x^*}\}=({}^*v)(x).
\end{align*}
Similarly, by part 3 of Proposition \ref{thm:duality} we have, for any $x\not\in C_u$, \[ -\infty=u(x)=\inf\sets{\UU_x^*(m)}{m\in\D}=({}^*v)(x). \]

Since $u$ and ${}^*v$ agree everywhere, except possibly on the boundary of $C_u$, it follows that $\cl u=\cl({}^*v)={}^*v$. Since $u$ is proper concave and $v$ is lower semi-continuous and proper convex, it follows that $u^*=(\cl(u))^*=({}^*v)^*=\cl(v)=v$ (c.f. Definition \ref{def:dualfunctionals}).\qed
\end{proof}

\begin{corollary}\label{thm:RAEimplAS}
Suppose that Assumption \ref{ass:proper_concave} holds. If either $U$ is bounded from above, or $U^*$ satisfies the growth condition \eqref{eqn:growth} then both $U$ and the value function $u$ are asymptotically satiable.
\end{corollary}

\begin{proof}
If $U$ is bounded from above then $U^*(0)=\sup_{x\in\RR_+^d}U(x)<\infty$, thus $0\in\dom(U^*)$ and hence $U$ is asymptotically satiable by Proposition \ref{thm:satiable}. Similarly, $u$ must also bounded from above in this case, and hence also asymptotically satiable.

Suppose that $U^*$ satisfies the growth condition. By Lemma \ref{thm:dual_functional} and the proper convexity of $U^*$, there exists an $x^*\in\interior(\RR_+^d)$ such that $U^*(x^*)<\infty$. It follows immediately from the growth condition that $\epsilon x^*\in\dom(U^*)$. Taking the limit as $\epsilon\rightarrow0$ shows that $0\in\cl(\dom(U^*))$, and hence $U$ is asymptotically satiable by Proposition \ref{thm:satiable}. We argue similarly to show that $u$ is asymptotically satiable. From part 2 of Proposition \ref{thm:duality}, and Lemma \ref{thm:repVV} we may choose any $m$ in the nonempty set $\D\cap\dom(\UU_0^*)\neq\emptyset$ (any minimizer in a dual problem with $x\in\interior(C_u)$ will do). Let $x^*:=m(\Omega)$, and let $\epsilon\in(0,1)$. Recall that $\proj:\RR^D\rightarrow\RR^d$ is defined by \eqref{eqn:def_proj}. By Proposition \ref{thm:corol}, Lemma \ref{thm:repVV} and \eqref{eqn:growth},
\begin{align*}
u^*(\epsilon x^*) &\le \UU_0^*(\epsilon m) = \E{\tilde U^*\left(\epsilon\radnik[m^c]\right)} = \E{U^*\left(\epsilon\proj\left(\radnik[m^c]\right)\right)} \\
&\le\zeta(\epsilon)\left(\E{U^*\left(\proj\left(\radnik[m^c]\right)\right)^+}+1\right) \\
&=\zeta(\epsilon)\left(\E{\tilde U^*\left(\radnik[m^c]\right)^+}+1\right)<\infty.
\end{align*}
We have shown that $\epsilon x^*\in\dom(u^*)$. Taking the limit as $\epsilon\rightarrow0$ shows that $0\in\cl(\dom(u^*))$, and hence $u$ is asymptotically satiable by Proposition \ref{thm:satiable}.

Note that if $U$ is bounded from above then $U^*$ satisfies the growth condition (see Remark \ref{rem:something}), and we could have used this to prove that $U$ and $u$ are asymptotically satiable. However, arguing this way would have been over-complicated.\qed
\end{proof}

Recall that if $x\in\interior(\dom(u))=\interior(C_u)$ then the existence of a minimizer $\hat m_x\in\D\cap\dom(\UU_x^*)$ in the dual problem \eqref{eqn:dual} is guaranteed by part 2 of Proposition \ref{thm:duality}. We now collect some of the properties of the minimizer. 

\begin{corollary}\label{thm:massofminimizer}
  Suppose that Assumptions \ref{ass:basic} and \ref{ass:proper_concave} hold. Given any $x\in\interior(\dom(u))$ and a minimizer $\hat m_x$ for the dual problem we have $\hat m_x(\Omega)\in\partial u(x)$.
\end{corollary}

\begin{proof}
  Define $x^*=\hat m_x(\Omega)$. Then by Proposition \ref{thm:duality}, Lemma \ref{thm:repVV} and Proposition \ref{thm:corol}
  \begin{align*}
    u(x) &= \UU_x^*(\hat m_x) = \UU_0^*(\hat m_x)+\ip x{x^*} \\
    &\ge \min\sets{\UU_0^*(m)}{m\in\D, m(\Omega)=x^*}+\ip x{x^*} \\
    &= u^*(x^*)+\ip x{x^*}.
  \end{align*}
  It now follows from \cite[Theorem 23.5]{Rock} that $x^*\in\partial u(x)$.\qed
\end{proof}

In the next result we shall see that (although the minimizer itself may not be unique) the first $d$ coordinate measures of the countably additive part of the minimizer are unique, and equivalent to $\PP$. The equivalence to $\PP$ is an essential ingredient in the paper, as it ensures that the random variable $\hat X_x$ in Proposition \ref{thm:variational} is well defined.
\begin{proposition}\label{thm:SCPSintval}
Suppose that Assumptions \ref{ass:basic} and \ref{ass:proper_concave} hold. Given any $x\in\interior(\dom(u))$, any minimizer $\hat m_x$ for the dual problem lies in the set $\P:=\sets{m\in\ba(\RR_+^D)}{\proj(\radnik[m^c])\text{ is }\interior(\RR_+^d)\text{-valued a.s.}}$, where $\proj:\RR^D\rightarrow\RR^d$ is defined by \eqref{eqn:def_proj}. Suppose that $\tilde m_x$ is another minimizer in the dual problem then $\proj(\radnik[\hat m_x^c])=\proj(\radnik[\tilde m_x^c])$ a.s. and $\hat m_x(x)=\tilde m_x(x)$.
\end{proposition}

\begin{remark}\label{thm:secondremark}
In the proofs of Proposition \ref{thm:SCPSintval} and Theorem \ref{thm:optimal} it will be useful to embed $\Z^s$ in $\D$ as follows. Given any $Z^s\in\Z^s$, we can construct a corresponding $m^s\in\ba(\RR_+^D)\cap\ca(\RR^D)$ by setting $m^s(A):=\E{Z_T^s\chi_A}$ for each $A\in\F_T$. It follows from \cite[Lemma 2.8]{CampScha} (which requires Assumption \ref{ass:SCPS}) that $m^s\in\D$. Note that $\radnik[m^s]=Z^s_T$ is $\interior(\RR_+^D)$-valued a.s. because $Z^s$ is a \emph{strictly} consistent price process.
\end{remark}

\begin{proof}[of Proposition \ref{thm:SCPSintval}]
Let $\partial\RR_+^d$ denote the boundary of $\RR_+^d$. Take $a\in\partial\RR_+^d$ and $b\in\interior(\RR_+^d)$. Recall from Lemma \ref{thm:essconc_esssmooth} and Proposition \ref{thm:satiable} that $U^*$ is strictly convex on $\interior(\RR_+^d)$, essentially smooth, and $\nabla U^*$ maps $\interior(\RR_+^d)$ into $-\interior(\RR_+^d)$. Since $U^*$ is essentially smooth, $|\nabla U^*(a+\lambda b)|\rightarrow\infty$ as $\lambda\rightarrow0$. Thus, by convexity of $U^*$,
\begin{equation}\label{eqn:useful}
\lim_{\lambda\searrow0}\frac{U^*(a+\lambda b)-U^*(a)}\lambda\le\lim_{\lambda\searrow0}\ip{\nabla U^*(a+\lambda b)}b=-\infty.
\end{equation}

From Lemma \ref{thm:repVV}, $\hat m_x\in\ba(\RR_+^D)$ and $\radnik[\hat m_x^c]$ is $\RR_+^D$-valued a.s. Suppose, for a contradiction, that $\hat m_x\not\in\P$. Then the event $A:=\big\{\proj(\radnik[\hat m_x^c])\in\partial\RR_+^d\big\}$ is non-null under $\PP$. Choose any $Z^s\in\Z^s$ (which is nonempty by Assumption \ref{ass:SCPS}), and let $m^s\in\D\cap\P$ be the corresponding Euclidean vector measure (see Remark \ref{thm:secondremark}). For $\lambda>0$, define $m_\lambda:=\hat m_x+\lambda m^s\in\D$ and $\nu_\lambda:=\tilde U^*\big(\radnik[m_\lambda^c]\big)$. Since, by Lemma \ref{thm:dual_functional}, $\UU_0^*$ is decreasing with respect to the preorder induced by $\ba(\RR_+^D)$, we see that $m_\lambda\in\dom(\UU_0^*)$. Since $\tilde U^*$ is convex, the integrable random variables $(\nu_\lambda-\nu_0)/\lambda$ are monotone increasing in $\lambda$. By the monotone convergence theorem and \eqref{eqn:useful}
\begin{align*}
&\lim_{\lambda\searrow0}\E{\chi_A^{}\left(\frac{\nu_\lambda-\nu_0}\lambda\right)}
= \E{\chi_A^{}\lim_{\lambda\searrow0}\left(\frac{\nu_\lambda-\nu_0}\lambda\right)} \\
&\qquad\qquad= \E{\chi_A^{}\lim_{\lambda\searrow0}\left(\frac{\tilde U^*(\radnik[\hat m_x^c]+\lambda \radnik[m^s])-\tilde U^*(\radnik[\hat m_x^c])}\lambda\right)} \\
&\qquad\qquad= \E{\chi_A^{}\lim_{\lambda\searrow0}\left(\frac{U^*\big(\proj(\radnik[\hat m_x^c])+\lambda\proj(\radnik[m^s])\big)-U^*\big(\proj(\radnik[\hat m_x^c])\big)}\lambda\right)}=-\infty.
\end{align*}
Hence $\lim_{\lambda\searrow0}\frac1\lambda\E{\nu_\lambda-\nu_0}=-\infty$. However, Lemma \ref{thm:repVV} and optimality of $\hat m_x$ imply that
\begin{align*}
  \E{\nu_\lambda-\nu_0}&=\E{\tilde U^*\left(\radnik[m_\lambda^c]\right)}-\E{\tilde U^*\left(\radnik[\hat m_x^c]\right)} \\
  &=\UU_x^*(m_\lambda)-m_\lambda(x)-\UU_x^*(\hat m_x)+\hat m_x(x)\ge-\lambda m^s(x).
\end{align*}
Therefore, for all $\lambda>0$, $\frac1\lambda\E{\nu_\lambda-\nu_0}\ge-m^s(x)$. This is the required contradiction.

Suppose for a contradiction that there exist solutions $\hat m_x,\tilde m_x$ to the dual problem such that\newline $\PP\big(\proj(\radnik[\hat m_x^c])\neq\proj(\radnik[\tilde m_x^c])\big)>0$. Defining $\bar m:=(\hat m_x+\tilde m_x)/2\in\D\cap\P$, strict convexity of $U^*$ on $\interior(\RR_+^d)$ implies that
\begin{align*}
&\E{\tilde U^*\left(\radnik[\bar m^c]\right)}+\bar m(x) = \E{U^*\left(\proj\bigg(\radnik[\bar m^c]\bigg)\right)}+\bar m(x) \\
&\qquad< \frac12\left\{\E{U^*\left(\proj\left(\radnik[\hat m_x^c]\right)\right)}+\hat m_x(x)\right\} \\
&\qquad\qquad\qquad\qquad\qquad\qquad\qquad\qquad+\frac12\left\{\E{U^*\left(\proj\left(\radnik[\tilde m_x^c]\right)\right)}+\tilde m_x(x)\right\}\\
&\qquad= \frac12\left\{\E{\tilde U^*\left(\radnik[\hat m_x^c]\right)}+\hat m_x(x)\right\}+\frac12\left\{\E{\tilde U^*\left(\radnik[\tilde m_x^c]\right)}+\tilde m_x(x)\right\}\\
&\qquad=\min_{m\in\D}\UU_x^*(m),
\end{align*}
which is the required contradiction. It follows immediately from Lemma \ref{thm:repVV} that $\hat m_x(x)=\tilde m_x(x)$.\qed
\end{proof}

\begin{proposition}[Variational Analysis]\label{thm:variational}
Suppose that Assumptions \ref{ass:basic}, \ref{ass:proper_concave} and \ref{ass:additional} hold. Given any $x\in\interior(\dom(u))$, let $\hat m_x\in\D\cap\dom(\UU_0^*)\cap\P$ denote an optimal dual measure, and define $\hat X_x:=\tilde I(\radnik[\hat m_x^c])$, where $\tilde I$ is defined by \eqref{def:functionI}. Then $\bigE{\bigip{\hat X_x}{\radnik[m^c]}}\le m(x)$ for all $m\in\D$, with equality for $m=\hat m_x$.
\end{proposition}

\begin{proof}
Take any $\tilde m\in\D\cap\dom(\UU_0^*)$. Since $\D$ and $\UU_0^*$ are convex, the measure $m_\lambda:=\lambda\tilde m+(1-\lambda)\hat m_x$ is again an element of $\D\cap\dom(\UU_0^*)$ for any $\lambda\in[0,1]$. The map $f:[0,1]\rightarrow\RR$ defined by $f(\lambda):=\UU_x^*(m_\lambda)$ is convex, and has a minimum at $0$. Therefore, by Lemma \ref{thm:repVV} and the Monotone Convergence Theorem,
\begin{align*}
0 &\le f'_+(0)=\lim_{\lambda\searrow0}\left\{\frac{f(\lambda)-f(0)}\lambda\right\} \\
&= \lim_{\lambda\searrow0}\left\{\E{\frac{\tilde U^*\left(\radnik[m_\lambda^c]\right)-\tilde U^*\left(\radnik[\hat m_x^c]\right)}\lambda} +\frac{m_\lambda(x)-\hat m_x(x)}\lambda\right\}\\
&= \E{\lim_{\lambda\searrow0}\left\{\frac{\tilde U^*\left(\radnik[m_\lambda^c]\right)-\tilde U^*\left(\radnik[\hat m_x^c]\right)}\lambda\right\}}+\tilde m(x)-\hat m_x(x)\\
&= \E{\ip{-\tilde I\left(\radnik[\hat m_x^c]\right)}{\radnik[\tilde m^c]-\radnik[\hat m_x^c]}}+\tilde m(x)-\hat m_x(x).
\end{align*}
Therefore
\begin{equation}\label{eqn:variational}
\E{\ip{\hat X_x}{\radnik[\tilde m^c]}}-\tilde m(x)\le\E{\ip{\hat X_x}{\radnik[\hat m_x^c]}}-\hat m_x(x).
\end{equation}
Assume now that $m\in\D$. It follows from Lemma \ref{thm:dual_functional} that $\UU_0^*$ is decreasing with respect to the preorder induced by $\ba(\RR_+^D)$, and hence $\tilde m:=\hat m_x+m\in\D\cap\dom(\UU_0^*)$. It follows from \eqref{eqn:variational} that
\begin{equation}\label{eqn:hhropgn}
\E{\ip{\hat X_x}{\radnik[m^c]}}\le m(x).
\end{equation}
By Proposition \ref{thm:satiable}, given any $\epsilon>0$ there exists an $x^*\in\dom(u^*)$ satisfying $\ip x{x^*}\le\epsilon$. Since $u^*(x^*)<\infty$, Proposition \ref{thm:corol} implies the existence of a $\tilde m\in\D\cap\dom(\UU_0^*)$ with $\tilde m(\Omega)=x^*$. By Lemma \ref{thm:repVV}, $\radnik[\tilde m^c]$ is $\RR_+^D$-valued a.s. Since $\hat X_x$ is also $\RR_+^D$-valued a.s. we have (using also \eqref{eqn:variational} and \eqref{eqn:hhropgn})
\begin{align*}
  -\epsilon\le-\ip x{x^*}=-\tilde m(x)&\le\E{\ip{\hat X_x}{\radnik[\tilde m^c]}}-\tilde m(x) \\
  &\le\E{\ip{\hat X_x}{\radnik[\hat m_x^c]}}-\hat m_x(x)\le0.
\end{align*}
Since $\epsilon>0$ is arbitrary, we have $\bigE{\bigip{\hat X_x}{\radnik[\hat m_x^c]}}=\hat m_x(x)$.\qed
\end{proof}

We now present our main theorem.
\begin{theorem}\label{thm:optimal}
Let $U:\RR^d\rightarrow[-\infty,\infty)$ be a utility function supported on $\RR_+^d$, which satisfies Assumption \ref{ass:basic}. Suppose in addition that Assumptions \ref{ass:proper_concave} and \ref{ass:additional} hold, and that the economy satisfies Assumption \ref{ass:SCPS}. Given any $x\in\interior(\dom(u))$, the optimal investment problem \eqref{eqn:primal} has a unique solution $\hat X_x:=\tilde I(\radnik[\hat m_x^c])$, where $\tilde I$ is defined by \eqref{def:functionI}, and where $\hat m_x$ is any dual optimizer from part 2 of Proposition \ref{thm:duality}.
\end{theorem}

\begin{proof}
Choose any $Z^s\in\Z^s$ (which is nonempty by Assumption \ref{ass:SCPS}), and let $m^s\in\D$ be the corresponding Euclidean vector measure (see Remark \ref{thm:secondremark}). It follows from Proposition \ref{thm:variational} that $\bigE{\bigip{\hat X_x}{Z^s_T}}=\bigE{\bigip{\hat X_x}{\radnik[m^s]}}\le m^s(x)=\ip x{Z^s_0}$. Theorem \ref{thm:superrep} implies that $\hat X_x\in\A{x}$. Furthermore, by Corollary \ref{thm:functionI}, Proposition \ref{thm:variational} and Lemma \ref{thm:repVV}, we have
\begin{align}\label{eqn:fenchequal}
\E{\tilde U(\hat X_x)}&=\E{\tilde U^*\left(\radnik[\hat m_x^c]\right)+\ip{\hat X_x}{\radnik[\hat m_x^c]}} = \E{\tilde U^*\left(\radnik[\hat m_x^c]\right)}+\hat m_x(x) \notag\\
&=\UU_x^*(\hat m_x).
\end{align}
It follows from part 1 of Proposition \ref{thm:duality} that $\hat X_x$ is an optimizer in the primal problem.

To show uniqueness, suppose for a contradiction that $\tilde X_x\in\A x$ is an optimizer in the primal problem such that $\PP(\tilde X_x\neq\hat X_x)>0$. Since $\tilde U$ has support cone $\RR_+^D$, $\tilde X_x$ must be $\RR_+^D$-valued a.s. By definition, $\hat X_x$ is $\interior(\RR_+^d)\times\RR_+^{D-d}$-valued a.s. We may assume without loss of generality that $\tilde X_x$ is also $\interior(\RR_+^d)\times\RR_+^{D-d}$-valued a.s., otherwise we can simply replace $\tilde X_x$ with the random variable $(\tilde X_x+\hat X_x)/2\in \A x$, which is $\interior(\RR_+^d)\times\RR_+^{D-d}$-valued a.s., and which is also an optimizer in the primal problem, due to concavity of $\tilde U$. Recall that $P:\RR^D\rightarrow\RR^d$ is defined by \eqref{eqn:def_proj}. There are two cases: Either (i) $\PP\big(\proj(\tilde X_x)\neq\proj(\hat X_x)\big)>0$ or (ii) $\PP\big(\ip{\smash{\tilde X_x}}{e^j}>0\big)>0$ for some $j\in\{d+1,\dots,D\}$.

(i) Define $\bar X:=(\tilde X_x+\hat X_x)/2\in\A x$. Since $U$ is strictly concave on $\interior(\RR_+^d)$,
\begin{align*}
  \E{\smash{\tilde U(\bar X_x)}}=\E{U(\proj(\bar X_x))}&>\frac12\{\E{\smash{U(\proj(\tilde X_x))}}+\E{\smash{U(\proj(\hat X_x))}}\} \\
&=\frac12\{\E{\smash{\tilde U(\tilde X_x)}}+\E{\smash{\tilde U(\hat X_x)}}\}=u(x),
\end{align*}
which is the required contradiction.

(ii) Let $j\in\set{d+1,\dots,D}$ be such that $\PP(\ip{\smash{\tilde X_x}}{e^j}>0)>0$. Define $\bar X_x:=\tilde X_x-Y$ where $Y:=\frac{\ip{\smash{\tilde X_x}}{e^j}}{\pi_T^{j,1}}(\pi_T^{j,1}e^j-e^1)$ is $K_T$-valued. Since $\ip{\smash{\bar X_x}}{e^1}\ge0$ a.s. and $\ip{\smash{\bar X_x}}{e^j}=0$ a.s., $\bar X_x$ is $\RR_+^D$-valued a.s. Hence $\bar X_x\in\A x$. Since $U$ is increasing with respect to $\presucc_{\RR_+^d}$ and strictly concave on $\interior(\RR_+^d)$, it must be strictly increasing on $\interior(\RR_+^d)$ with respect to $\presucc_{\RR_+^d}$. Hence
\begin{align*}
\E{\tilde U(\bar X_x)} &= \E{U\big(\proj(\tilde X_x)-\proj(Y)\big)} = \E{U\left(\proj(\tilde X_x)+\frac{\ip{\smash{\tilde X_x}}{e^j}}{\pi_T^{j,1}}e^1\right)} \\
&> \E{U\big(\proj(\tilde X_x)\big)} = \E{\tilde U(\tilde X_x)} = u(x),
\end{align*}
which is the required contradiction.\qed
\end{proof}

We finish this section by giving an example where the singular part, $\hat m_x^p$, of the dual minimizer is non-zero.
\begin{example}\label{exa:singularcomponent}
Let $S:=(S_0,S_1)$ be as defined in \cite[Example 5.1']{KramScha99}. That is, $S_0\equiv1$ and $S_1$ takes the values $(s_n)_{n=0}^\infty$ with probabilities $(p_n)_{n=0}^\infty$, where $s_0=2$, $s_n=1/n$ for $n\ge1$, $p_0=1-\alpha$ and $p_n=\alpha2^{-n}$, with $\alpha$ sufficiently small. 
This example can be modified to include frictions as follows: With $D=2$, we define the bid-ask process
\[ \Pi_0:=\begin{pmatrix}1 & S_0 \\ 2/S_0 & 1 \end{pmatrix}=\begin{pmatrix}1 & 1 \\ 2 & 1 \end{pmatrix}\qquad\text{and}\qquad\Pi_1:=\begin{pmatrix}1 & 2S_1 \\ 1/S_1 & 1 \end{pmatrix}, \]
and let $\A0$ denote the corresponding cone of admissible terminal portfolios with zero initial portfolio.

Note that under this model the $\RR^2$-valued price process $(1,S_t)$, $t=0,1$, is now a shadow price for the bond and stock. In relation to this shadow price process, at time $t=0$, trading from the bond to the stock is frictionless, while trading in the opposite direction incurs costs. At time $t=1$, however, trading from the stock to the bond is now frictionless, while trading from bond to stock incurs costs. 

With $d=1$, we set $U(x):=\ln x$. We define the extended utility function $\tilde U:\RR^2\rightarrow[-\infty,\infty)$ by \eqref{eqn:extension}, and the value function $u:\RR^2\rightarrow[-\infty,\infty)$ by \eqref{eqn:primal}. Since $1=d<D=2$, the extended utility function effectively forces the investor to close out their position in the stock at maturity, in order to derive the maximum possible utility from their terminal portfolio.

Suppose we are given an initial portfolio $x=(x_0,x_1)$. In the frictionless case, since $S_0\equiv1$ we may immediately trade at time $0$ to the portfolio $(x_0+x_1,0)$, and hence the maximum expected utility is given by $\tilde u(x):=u^{\text{KS}}(x_1+x_2)$, where $u^{\text{KS}}(x):=\ln x+\E{\ln S_1}$ is the value function obtained in \cite[Example 5.1']{KramScha99}. However, if we introduce frictions as described above, this only serves to decrease the terminal wealth, and hence the associated utility. Thus $u(x)\le \tilde u(x)$.

We shall now see that $u$ and $\tilde u$ are equal whenever $x_1>0$ and $x_2\ge-x_1$. We claim that $X_x:=((x_1+x_2)S_1,0)\in\A{x}$. Indeed, to reach this terminal portfolio from the initial portfolio $x=(x_1,x_2)$, one can trade to $(0,x_1+x_2)$ at time $0$ and then at time $1$, $X_x$ can be reached by liquidating to the bond. Thus
\begin{align*}
  \E{\smash{\tilde U(X_x)}}&=\E{U((x_1+x_2)S_1)}=\ln(x_1+x_2)+\E{\ln S_1}=u^{\text{KS}}(x_1+x_2) \\
  &= \tilde u(x)\ge u(x).
\end{align*}
Hence $X_x=\hat X_x$ is optimal and $u(x)=\ln(x_1+x_2)+\E{\ln S_1}$.

Now fix $x=(1,0)$. Let $\hat m=\hat m_x$ denote the minimizer in the dual problem. By Corollary \ref{thm:massofminimizer}, $\hat m(\Omega)\in\partial u(1,0)=\set{(1,1)}$. In particular the first coordinate measure $\hat m_1:=\ip{e^1}{\hat m}$ satsifies $\hat m_1(\Omega)=1$. By Theorem \ref{thm:optimal},
\[ (S_1,0)=X_{(1,0)}=\hat X_{(1,0)}=\left(-(U^*)'\left(\radnik[\hat m_1^c]\right),0\right). \]
Hence
\[ \radnik[\hat m_1^c]=U'(S_1)=\frac1{S_1}. \]
Referring back to \cite[Example 5.1']{KramScha99}, we see that
\[ \hat m_1^c(\Omega)=\E{\radnik[\hat m_1^c]}=\E{\frac1{S_1}}<1. \]
Since $\hat m_1(\Omega)=1$ and $\hat m_1^c(\Omega)<1$, it be the case that $\hat m_1^p(\Omega)\neq0$.\qed
\end{example}

\section{The liquidation case}\label{sec:liquidation}

In many papers dealing with optimal investment under transaction costs, it is assumed that the agent liquidates their assets at the close of trading to a given reference asset, which is chosen as a num\'eraire at time $t=0$. The reader is referred especially to \cite{Kaba}, \cite{DeelPhamTouz}, \cite{Bouc} and the references therein. In this subsection, we show that our optimal investment problem is equivalent to maximizing expected utility from liquidation of the terminal portfolio, thus avoiding the delicate issue of using a non-smooth utility function as in \cite{DeelPhamTouz}.

\begin{definition}\label{def:luciano}
Let $U$ be a utility function supported on $\RR_+^d$ (see Definition \ref{def:utility}) which satisfies Assumption \ref{ass:basic}. The \emph{terminal liquidation utility functional} corresponding to $U$ is defined\footnote{Clearly, the set over which we are optimizing in \eqref{eqn:indU} is a.s. nonempty (the zero vector belongs to it) and compact in $\RR^d _+$. Since $U$ is upper semi-continuous, this justifies the use of the maximum for almost every $\omega$.}
by
\begin{equation}\label{eqn:indU}
 \bar U(W) := \max\sets{U(\xi)}{\xi\in \RR_+ ^d ,\;(\xi,\underline0)-W\in -K_T },\quad W\in L^0(K_T,\F_{T-}),
\end{equation}
where $\underline 0$ denotes the zero vector in $\RR^{D-d}$.
\end{definition}

Given $W\in L^0(K_T,\F_{T-})$, the random quantity $\bar U(W)$ models the best an agent can do if, at time $T$, they decide to liquidate their portfolio at time $T-$ to the $d$ consumption goods according to the terminal solvency cone $K_T$. Observe that it is natural to consider only those random variables $W$ that belong to $K_T$ a.s., since $W$ represents agent's portfolio at time $T-$ resulting from an admissible portfolio $V \in \mathcal A^x$ for some initial endowment $x$. Indeed, $V_{T-}=(V_{T-}-V_T) + V_T$ where $V_{T-}-V_T \in K_T$ and, without loss of generality, $V_T \in \RR_+ ^D$, so that $V_{T-}$ belongs a.s. to $K_T + \RR_+ ^D = K_T$.

\begin{remark}\label{rem:meas-selection}
Before stating the main results of this section, we notice that for any $W\in L^0(K_T,\F_{T-})$ the liquidation functional $\bar U(W)$ defined by \eqref{eqn:indU} admits a \emph{measurable} maximum $\hat\xi$ (i.e. the set of maximizers admits a measurable selector). To prove this, note that we can reformulate the terminal liquidation functional $\bar U(W)$ as
\[ m(\omega):=\max\sets{f(\omega,\xi)}{\xi\in\phi(\omega)}, \]
where $f:\Omega\times \RR_+ ^d \rightarrow\RR$ is defined by $f(\omega,\xi):=U(\xi)$, and $\phi:\Omega\twoheadrightarrow\RR_+ ^d$ is defined by $\phi(\omega):=\sets{\xi\in\dom(U)}{(\xi,\underline0)-W(\omega)\in-K_T(\omega)}$. Since $W\in K_T$ a.s., $\phi$ has nonempty and compact values a.s. It follows from \cite[Lemmas 18.3 and 18.7]{AlipBord} that $\phi$ is weakly measurable. Since $U$ is upper semi-continuous, $f$ is Carath\'eodory. Thus $\phi$ and $f$ satisfy the conditions of the measurable maximum theorem \cite[Theorem 18.19]{AlipBord} except from the fact that $f$ can take the value $-\infty$. Nonetheless \cite[Theorem 18.9]{AlipBord} can be applied\footnote{For the sake of clarity, we notice that even though \cite[Theorem 18.19]{AlipBord} is stated only for finite-valued functions $f$, it can be applied to functions taking possibly the value $-\infty$ as follows: Let $\psi$ be an order-preserving homeomorphism mapping $[-\infty,\infty)$ into $[0,1)$. One can apply \cite[Theorem 18.19]{AlipBord} to the function $\psi\circ f$ to get a measurable maximizer. Since $\psi$ is order-preserving, such a maximizer coincides with that of our original maximization problem.}  
so that, in particular, the argmax correspondence of maximizers $\mu:\Omega\twoheadrightarrow\RR_+ ^d$ defined by $\mu(\omega):=\sets{\xi\in\phi(\omega)}{f(\omega,\xi)=m(\omega)}$ admits a measurable selector $\hat\xi:\Omega\rightarrow\RR_+ ^d$.
\end{remark}

The following propositions are the two main results of this section: In Proposition \ref{thm:liquid-equiv} we show that the value function of the original problem coincides with the supremum of the expected liquidation utility functional. In Proposition \ref{thm:liquid-attained} we go on to show that both problems essentially have a common optimizer.

\begin{proposition}\label{thm:liquid-equiv}
Let $x\in\RR^D$ be a given initial endowment. Then
\begin{equation}\label{eqn:dirUindU}
u(x)=\sup_{W\in\mathcal A_{T-}^x}\E{\bar U(W)},
\end{equation}
where $\mathcal A_{T-}^x:=\sets{V_{T-}}{V\in\mathcal A^x}$.
\end{proposition}

\begin{proof}
  First, we prove inequality `$\le$'. Let $V$ be a given admissible portfolio process such that $V_0=x$. We assume without loss of generality that $V_T\in\RR_+^D$ a.s. It follows from \cite[Lemma 2.8]{CampScha} and Assumption \ref{ass:contPi} that $(\proj(V_T),\underline0)-V_{T-}\in-K_T$ a.s., where $\proj:\RR^D\rightarrow\RR^d$ is defined by \eqref{eqn:def_proj}. Hence, by definition of $\tilde U$ and $\bar U$, we have
\[ \tilde U(V_T) = U(\proj(V_T)) \le \sup\sets{U(\xi )}{\xi\in \RR^d,\;(\xi,\underline 0)-V_{T-}\in-K_T } = \bar U(V_{T-}). \]
Hence the desired inequality follows.

For the opposite inequality `$\ge$', let $V\in\mathcal A^x$. By Remark \ref{rem:meas-selection} there exists a $\F_T$-measurable solution $\hat\xi$ to the optimization problem \eqref{eqn:indU} when $W=V_{T-}$. Indeed, as we have already noticed, $V_{T-}$ belongs to $K_T$ and thus the maximizer $\hat \xi$ is well-defined. Moreover, the strict concavity of $U$ implies that such a maximizer is a.s. unique.

We claim that $(\hat\xi,\underline 0)$ belongs to $\A x$. Indeed, $(\hat\xi,\underline 0)$ is the terminal value of the portfolio process $V'$ defined as $V'_t=V_t+((\hat\xi,\underline 0)-V_T)\chi_{\set{t=T}}$, which clearly belongs to $\mathcal A^x$ because over $[0,T)$ it coincides with $V$ which is admissible and at $T$ the last trade equals $\Delta V^\prime _T =V'_T-V'_{T-}=(\hat\xi,\underline 0)-V_{T-}\in-K_T$ a.s. As a consequence, one has
\[ u(x)\ge\E{ U(\hat\xi )}=\E{\bar U(V_{T-})} \]
which gives the result.\qed
\end{proof}

\begin{proposition}\label{thm:liquid-attained}
The supremum in \eqref{eqn:dirUindU} is attained. Moreover, given any maximizer $\hat W$ in \eqref{eqn:dirUindU}, let $\hat\xi=\hat\xi(\hat W)$ be any maximizer in the optimization problem $\bar U(\hat W)$ and let $\hat X_x$ be the unique maximizer in the primal problem \eqref{eqn:primal}. Then $(\hat\xi(\hat W),\underline 0)=\hat X_x$ a.s.
\end{proposition}

\begin{proof}
Since $\hat X_x\in\A x$, there exists an admissible $V$ such that $V_0=x$ and $V_T=\hat X_x$. Define $\hat W:=V_{T-}$, and $\hat\xi:=P(\hat X_x)$. By \cite[Lemma 2.8]{CampScha}, $(\hat\xi,\underline 0)-\hat W=\hat X_x-V_{T-}=V_T-V_{T-}\in-K_T$ a.s. Now
\[ \E{\bar U(\hat W)}\ge\E{U(\hat\xi)}=\E{\tilde U(\hat X_x)}=u(x). \]
Therefore by Proposition \ref{thm:liquid-equiv}, $\hat W$ is optimal in \eqref{eqn:dirUindU}. Now suppose that $\tilde W$ is any maximizer in \eqref{eqn:dirUindU}, and let $\tilde\xi=\tilde\xi(\tilde W)$ be the corresponding maximizer in $\bar U(\tilde W)$. Define $\tilde X_x:=(\tilde\xi,\underline 0)\in\A x$. Then
\[ \E{\tilde U(\tilde X_x)}=\E{U(\tilde\xi)}=\E{\bar U(\tilde W)}=u(x). \]
By Theorem \ref{thm:optimal}, $(\tilde\xi(\tilde W),\underline0)=\tilde X_x=\hat X_x$ a.s.\qed
\end{proof}

\begin{example}[Liquidation to the first asset] Take $d=1$, i.e. at the end the agent is interested in consuming only the first good. In this case a direct computation leads to the following expression for $\bar U$:
\[ \bar U(W)=U(l(W)), \]
where $l$ is the liquidation functional expressed in \emph{physical units}, defined as follows
\begin{equation}\label{liquid}
l(W)=\sup\sets{\xi\in \RR_+ }{(\xi,\underline 0)-W\in -K_T },\quad W\in L^0(K_T,\F_{T-}).
\end{equation}
Observe that while $U$ is smooth, the corresponding indirect utility function $\bar U$ need not be. The previous proposition can be rewritten as
\[ u(x) = \sup_{W\in\mathcal A_{T-}^x}\E{U(l(W))}. \]

We note that the function $l$ given in (\ref{liquid}) is the analogue (in our framework) of the liquidation function as defined, e.g., in the papers \cite{DeelPhamTouz} and \cite{Bouc}, where all quantities are expressed in terms of a fixed num\'eraire.\end{example} 

\section{Appendix}

\subsection{Lagrange duality}

The Lagrange duality theorem is the central ingredient in the proof of Proposition \ref{thm:duality}. Part 1 of the theorem below is essentially a reformulation of \cite[Theorem 8.6.1]{Luen} in terms of concave functionals which may take the value $-\infty$, as opposed to real-valued convex functionals. We have also added part 2 to cover the case where the optimization is degenerate.
\begin{theorem}[Lagrange duality theorem]\label{thm:lagrange}
Let $\X$ denote a normed\footnote{It is worth noting that the Lagrange duality theorem is also true if $\X$ is simply a topological vector space. We do not need the strengthened version of the result however, so we restrict ourselves to the case where $\X$ is a normed vector space.} 
vector space, let $\C$ be a nonempty convex cone in $\X$, let $\D:=(-\C)^*$, and let $\U:\X\rightarrow[-\infty,\infty)$ be a proper concave functional.
\begin{enumerate}
\item Suppose there exists a $p\in\interior(\C)$ and an $x\in\C$ such that $\U(x+p)>-\infty$, and $\sup_{x\in\C}\U(x) < \infty$. Then
    \[ \sup_{x\in\C}\U(x)=\min_{x^*\in\D}\U^*(x^*)\in\RR. \]
\item Suppose there exists a $p\in\interior(\C)$ such that $\U(x-p)=-\infty$ for all $x\in\C$ and there exists $x_1^*\in\D$ such that $\U^*(x_1^*)<\infty$. Then
    \[ \sup_{x\in\C}\U(x)=\inf_{x^*\in\D}\U^*(x^*)=-\infty. \]
\end{enumerate}
\end{theorem}

\begin{proof}
Note first that for any $x^*\in\D$ we have
\[ \sup_{x\in\C}\U(x)\le\sup_{x\in\C}\{\U(x)-\ip x{x^*}\}\le\sup_{x\in \X}\{\U(x)-\ip x{x^*}\}=\U^*(x^*). \]
\begin{enumerate}
  \item Following the notation of \cite[\S8]{Luen}, we set $X=Z=\X$, $\Omega=\dom(\U)$, and let $G:X\rightarrow Z$ be the identity operator. Let $P=-\C$ be the positive cone of $Z$, so that the dual, positive cone of $Z^*$ is $\D$. By the hypothesis of part 1, the point $x_1:=x+p$ lies both in the effective domain of $\U$ and in the interior of $\C$; in the notation of \cite[\S8]{Luen}, $x_1\in\Omega$ satisfies $G(x_1)<\theta$. Let $f$ be the restriction of $-\U$ to $\Omega$, thus $f$ is a real-valued convex functional defined on the convex subset $\Omega$ of $X$. It is easy to verify that the concave dual of $f$ is $\phi=-\U^*$. Applying \cite[Theorem 8.6.1]{Luen} gives
\begin{align*}
  \sup_{x\in\C}\U(x)&=-\inf\sets{f(x)}{G(x)\le\theta,\;x\in\Omega} \\
&=-\max\sets{\phi(x^*)}{x^*\ge\theta}=\min_{x^*\in\D}\U^*(x^*)\in\RR.
\end{align*}

  \item First note that
\[ \sup_{x\in\C}\U(x)\le\sup_{x\in-p+\C}\U(x)=\sup_{x\in\C}\U(x-p)=-\infty. \]
Furthermore, by the hypothesis of part 2, $\C$ and $S:=\sets{x'\in \X}{\U(x'-p)>-\infty}$ are disjoint, nonempty, convex sets. Since $\C$ is a convex cone which contains an interior point, \cite[Theorem V.2.8]{DunfSchw} implies the existence of a non-zero $x_0^*\in \X^*$ such that
\begin{equation}\label{eqn:working}
\ip x{x_0^*}\le0\le\ip{x'}{x_0^*}
\end{equation}
for all $x\in\C$ and all $x'\in S$. This implies that $x_0^*\in\D$.

Note that since $x_0^*\in\D$ and $p\in\C$, we have $\ip p{x_0^*}\le0$. We claim that $\ip p{x_0^*}<0$. Indeed, suppose for a contradiction that $\ip p{x_0^*}=0$. Since $x_0^*\neq0$, there exists an $x'\in \X$ such that $\ip{x'}{x_0^*}>0$. Since $p$ is an interior point of $\C$, by continuity of scalar multiplication there exists an $\epsilon>0$ such that $x'':=p+\epsilon x'\in\C$. Therefore $\ip{x''}{x_0^*}=\epsilon\ip{x'}{x_0^*}>0$, which contradicts the fact that $x_0^*\in\D$.

Given any $x\in\dom \U$, we have $x':=p+x\in S$. Hence by \eqref{eqn:working} we have
\begin{equation}\label{eqn:cough}
-\ip x{x_0^*}=\ip p{x_0^*}-\ip{x'}{x_0^*}\le\ip p{x_0^*}.
\end{equation}
Given any $\lambda>0$, note that $x_1^*+\lambda x_0^*\in\D$. It follows from the definition of $\U^*$ and \eqref{eqn:cough} that
\[ \U^*(x_1^*+\lambda x_0^*) = \sup_{x\in\dom(\U)}\{\U(x)-\ip x{x_1^*}-\lambda\ip x{x_0^*}\} \\
\le \U^*(x_1^*)+\lambda\ip p{x_0^*}. \]
Since $\U^*(x_1^*)<\infty$ and $\ip p{x_0^*}<0$ we may make the right-hand side arbitrarily negative by choosing $\lambda$ arbitrarily large. Therefore $\inf_{x^*\in\D}\U^*(x^*)=-\infty$.\qed
\end{enumerate}
\end{proof}

\subsection{Proofs of Auxiliary Results from Section \ref{sec:prelim}}

\subsubsection*{Proof of Lemma \ref{thm:nec_suf_ass_sat}}

Take any $\epsilon>0$ and suppose that there exists an $x\in\interior(\dom(U))$ such that $\partial U(x)\cap[0,\epsilon)^d\neq\emptyset$. By \cite[Corollary 23.5.2]{Rock}, $\partial(\cl(U))(x)=\partial U(x)$, and hence $U$ is asymptotically satiable.

Conversely, suppose that $U$ is essentially smooth and asymptotically satiable. By \cite[Theorem 7.4]{Rock}, $\cl(U)$ agrees with $U$ except perhaps at boundary points of $\dom(U)$. Therefore $\cl(U)$ is essentially smooth. Since $U$ is asymptotically satiable, given any $\epsilon>0$ there exists an $x\in\RR^d$ such that $\partial(\cl(U))(x)\cap[0,\epsilon)^d\neq\emptyset$. By \cite[Theorem 26.1]{Rock} we must have $x\in\interior(\dom(\cl(U)))=\interior(\dom(U))$, and $\nabla U(x)=\nabla(\cl(U))\in[0,\epsilon)^d$.\qed

\subsubsection*{Proof of Proposition \ref{thm:satiable}}

1 $\Rightarrow$ 2. For each $n\in\NN$ there exists an $x_n\in\RR^d$ such that $\partial(\cl(U))(x_n)\cap [0,1/n)^d\neq\emptyset$. Choose any $x_n^*\in\partial(\cl(U))(x_n)\cap [0,1/n)^d$. By \cite[Theorem 12.2 and Corollary 23.5.1]{Rock} we have $-x_n\in\partial (\cl(U)^*)(x_n^*)=\partial U^*(x_n^*)$ and hence, by \cite[Theorem 23.4]{Rock}, $x_n^*\in\dom(U^*)$. Since the sequence $(x_n^*)_{n\in\NN}$ converges to $0$, we have $0\in\cl(\dom(U^*))$.\vspace{0.5\baselineskip}

\noindent 2 $\Rightarrow$ 3. There exists a sequence $(x_n^*)_{n\in\NN}\subseteq\dom(U^*)$ such that $x_n^*\rightarrow0$ as $n\rightarrow\infty$. By Lemma, \ref{thm:dual_functional} $\dom(U^*)\subseteq (C_U)^*$. Take any $x^*\in\ri((C_U)^*)$. Since $x_n^*\rightarrow0$ as $n\rightarrow\infty$, the sequence $(x^*-x_n^*)_{n\in\NN}\subseteq\aff((C_U)^*)$ is eventually in $\ri((C_U)^*)$. Therefore $x^*\presucc_{(C_U)^*}x_n^*$ eventually, and since, by Lemma \ref{thm:dual_functional}, $U^*$ is decreasing with respect to $\presucc_{(C_U)^*}$, this implies that $x^*\in\dom(U^*)$. We have therefore shown that $\ri((C_U)^*)\subseteq\dom(U^*)$. By \cite[Corollary 6.3.1]{Rock}, this, together with the fact that $\dom(U^*)\subseteq (C_U)^*$, shows that $\cl(\dom(U^*))=(C_U)^*$.\vspace{0.5\baselineskip}

\noindent 3 $\Rightarrow$ 4. Obvious.\vspace{0.5\baselineskip}

\noindent 4 $\Rightarrow$ 1. By \cite[Corollary 6.3.1]{Rock}, $\cl(\dom(U^*))=\cl(\ri(\dom(U^*)))$. Since $\cl(\dom(U^*))$ is a convex cone, given any $\epsilon>0$ we may find a $x^*\in\ri(\dom(U^*))\cap [0,\epsilon)^d$. By \cite[Theorem 23.4]{Rock}, $\partial U^*(x^*)\neq\emptyset$. Choose any $x\in-\partial U^*(x^*)$. By \cite[Theorem 12.2 and Corollary 23.5.1]{Rock}, $x^*\in\partial(\cl(U))(x)$. Since $x^*\in\partial(\cl(U))(x)\cap [0,\epsilon)^d$, we have shown 1.\qed

\subsubsection*{Proof of Proposition \ref{thm:oldgrowth}}

Since $U$ satisfies \eqref{eqn:oldasymptelast} there exist $\beta>0$, $c\in\RR$, and $r>0$ such that for all $x\in\interior(\RR_+^d)$ satisfying $|x|\ge r$ we have $U(x)\ge(1+1/\beta)\ip x{\nabla U(x)}-c$. Let $\one\in\RR^d$ denote the vector whose entries are all equal to $1$. Define $x_r:=r\one$, and $x_r^*:=\nabla U(x_r)$.

Take any $x^*\in\interior(\RR_+^d)$ and $\epsilon\in(0,1]$. We consider two cases, (i) $x^*\presucc_{\RR_+^d}x_r^*$ and (ii) $x^*\not\presucc_{\RR_+^d}x_r^*$.
\begin{itemize}
\item[(i)] In this case $\epsilon x^*\presucc_{\RR_+^d}\epsilon x_r^*$, so by Lemma \ref{thm:dual_functional}, $U^*(\epsilon x^*)\le U^*(\epsilon x_r^*)$.
\item[(ii)] Since $U$ is asymptotically satiable, Proposition \ref{thm:satiable} shows that $\epsilon x^*\in\interior(\dom(U^*))$. By Lemma \ref{thm:essconc_esssmooth} we may define $x_\epsilon:=-\nabla U^*(\epsilon x^*)$. We claim that $|x_\epsilon|\ge r$. Indeed, suppose for a contradiction that $|x_\epsilon|<r$. Then $x_\epsilon\prepreq_{\RR_+^d} x_r$, so by Lemmas \ref{thm:essconc_esssmooth} and \ref{thm:mvra}, $x^*\presucc_{\RR_+^d}\epsilon x^*=\nabla U(x_\epsilon)\presucc_{\RR_+^d}\nabla U(x_r)=x_r^*$, which is the required contradiction. Therefore, by Corollary \ref{thm:functionI},
\begin{align}\label{eqn:oldaltgrowth}
U^*(\epsilon x^*) &= U(x_\epsilon)-\ip {x_\epsilon}{\epsilon x^*} \notag\\
&\ge (1+1/\beta)\ip {x_\epsilon}{\nabla U(x_\epsilon)}-c-\ip {x_\epsilon}{\epsilon x^*} \notag\\
&= -\frac1\beta\ip{\nabla U^*(\epsilon x^*)}{\epsilon x^*}-c.
\end{align}
Define the function $F:(0,1]\rightarrow\RR$ by $F(\epsilon):=\epsilon^\beta(U^*(\epsilon x^*)+c)$. Using \eqref{eqn:oldaltgrowth}, we see that
\[ F'(\epsilon)=\beta\epsilon^{\beta-1}(U^*(\epsilon x^*)+c+\ip{\nabla U^*(\epsilon x^*)}{\epsilon x^*}/\beta)\ge0. \]
Hence $U^*(\epsilon x^*)=\epsilon^{-\beta}F(\epsilon)-c\le\epsilon^{-\beta}F(1)-c=\epsilon^{-\beta}U^*(x^*)+(\epsilon^{-\beta}-1)c$.
\end{itemize}
The result follows by setting $\zeta(\epsilon):=\max\set{\epsilon^{-\beta},(\epsilon^{-\beta}-1)c,U^*(\epsilon x_r^*),0}$.\qed

\subsubsection*{Proof of Lemma \ref{thm:addUgrowth}}
Applying Proposition \ref{thm:oldgrowth} with $d=1$, for each $i\in\set{1,\dots,d}$ there exists a function $\zeta_i:(0,1]\rightarrow(0,\infty)$ such that for all $\epsilon\in(0,1]$ and all $x_i^*>0$
\[ U_i^*(\epsilon x_i^*)\le\zeta_i(\epsilon)(U_i(x_i^*)^++1). \]
It follows that for $x^*\in\interior(\RR_+^d)$,
\[ U^*(\epsilon x^*)=\sum_{i=1}^d U_i^*(\epsilon x_i^*)\le\sum_{i=1}^d\zeta_i(\epsilon)(U_i^*(x_i^*)^++1)\le\max_{i=1,\dots,d}\zeta_i(\epsilon)\left(\sum_{i=1}^d U_i^*(x_i^*)^++d\right). \]
Since $\inf\sets{U(x)}{x\in\interior(\RR_+^d)}>-\infty$, it follows that $a_i:=\inf\sets{U_i(x_i)}{x_i\in\interior(\RR_+)}>-\infty$ for each $i$. Moreover, since $U_i^*(x_i^*)^+=U_i^*(x_i^*)+U_i^*(x_i^*)^-\le U_i^*(x_i^*)+a_i^-$ we have
\[ \sum_{i=1}^d U_i^*(x_i^*)^+ \le U^*(x^*)+\sum_{i=1}^da_i^-\le U^*(x^*)^++\sum_{i=1}^da_i^-. \]
The growth condition follows by setting $\zeta(\epsilon)=\max_{i=1,\dots,d}\zeta_i(\epsilon)(\sum_{i=1}^da_i^-+d)$.\qed

\subsubsection*{Acknowledgement}

The authors thank Alfred M\"uller, Dmitry Kramkov, Paolo Guasoni, Steve Shreve, Beate Zimmer, Joe Diestel, John Wright and Jim Brooks for discussions about various topics relating to the paper. We would also like to thank an anonymous referee for comments which improved the presentation of the paper.

The first named author thanks the ``Chair Les Particuliers Face aux Risques'', Fondation du Risque (Groupama-ENSAE-Dauphine), the GIP-ANR ``Croyances'' project and the ``Chair Finance and Sustainable Development'' sponsored by EDF for their support. The second named author gratefully acknowledges partial support from EPSRC grant GR/S80202.

\end{document}